\documentclass[12pt, twoside, leqno]{article}

\usepackage{amsmath,amsthm}
\usepackage{amssymb}

\usepackage{enumitem}

\usepackage{graphicx}

\usepackage[T1]{fontenc}

\newtheorem{theorem}{Theorem}[section]
\newtheorem{corollary}[theorem]{Corollary}
\newtheorem{lemma}[theorem]{Lemma}
\newtheorem{proposition}[theorem]{Proposition}

\theoremstyle{definition}
\newtheorem{definition}[theorem]{Definition}
\newtheorem{remark}[theorem]{Remark}

\numberwithin{equation}{section}

\newcommand{\R}{\mathbb R}
\newcommand{\N}{\mathbb N}

\renewcommand{\u}{u^{\epsilon,\delta}}
\renewcommand{\div}{\rm{div}}
\usepackage{tikz}
\usepackage{amsmath}

\newcommand{\E}{\mathbb E}

\begin{document}

\title{Weak solutions to granular medium equations via stochastic particle approximations}

%%% Give just the affiliation and not the detailed postal address.

\author{Łukasz Wodnicki\\
Institute of Fundamental
Technological Research\\ 
Polish Academy of Sciences\\
 Poland\\and\\
Center of Trustworthy AI for Life Sciences\\
University
of Warsaw\\
 Poland\\
E-mail: 
l.wodnicki@uw.edu.pl
}
\date{}
\maketitle

%% Classification and key words; note that the 2020 classification is used.
%% Only one classification should be declared as primary.

\renewcommand{\thefootnote}{}

\footnote{2026 \emph{Mathematics Subject Classification}: Primary 60H10; Secondary 60K35, 35K65, 82C31.}

\footnote{\emph{Key words and phrases}: Particle approximation, stochastic differential equations.}

\renewcommand{\thefootnote}{\arabic{footnote}}
\setcounter{footnote}{0}

%%%%%%%%

\begin{abstract}
We establish weak existence for McKean--Vlasov stochastic differential
equations associated with granular medium equations with porous medium
diffusion, for $m>1$ in arbitrary dimension. The result covers a class
of nonconvex potentials with locally Lipschitz gradients of polynomial
growth. Our particle approximation generates nonlinear diffusion
through density-dependent noise and allows initial densities with
finite Boltzmann entropy but infinite internal energy.
\end{abstract}

\section{Introduction}

The granular medium equation with porous medium diffusion has been
studied for many years. However, a detailed proof of existence of weak
solutions to the associated McKean--Vlasov stochastic differential
equation does not appear to be readily available in the literature.
The aim of this work is to fill this gap. Specifically, we establish
the existence of weak solutions to a stochastic differential equation
(SDE) whose time-marginal densities solve the following partial
differential equation (PDE):
\begin{align}\label{granular medium equation}
    \left\{
    \begin{array}{ll}
        \partial_t u
        = \div(u\nabla V)
        + \div\bigl(u(\nabla W*u)\bigr)
        + \Delta u^m
        & \text{in }\R^d\times(0,T),\\
        u(\cdot,0)=u_0
        & \text{in }\R^d.
    \end{array}
    \right.
    \tag{PDE}
\end{align}
Here, $m>1$, $V$ and $W$ are the confinement and interaction potentials,
respectively, and $*$ denotes convolution in the spatial variable.
For each $t\in[0,T)$, the function $u_t=u(t,\cdot)$ is a probability
density on $\R^d$.

More precisely, we construct a process $(Y_t)_{t\in[0,T)}$ that solves,
in the weak sense, the McKean--Vlasov SDE
\begin{align}\label{SDE}
    \left\{
    \begin{array}{ll}
        dY_t
        = -\bigl(\nabla V(Y_t)+(\nabla W*u_t)(Y_t)\bigr)\,dt
        + \sqrt{2}\,u_t(Y_t)^{\frac{m-1}{2}}\,dB_t,
        & t\in(0,T),\\
        \operatorname{Law}(Y_t)(dx)=u_t(x)\,dx,\\
        Y_0\sim u_0(x)\,dx,
    \end{array}
    \right.
    \tag{SDE}
\end{align}
where $(B_t)_{t\in[0,T)}$ is a standard $d$-dimensional Brownian motion.

In the one-dimensional case, equation \ref{granular medium equation} was introduced to model granular media in \cite{B1,B2}. Consider many infinitesimal particles colliding inelastically, then $u(t,x)$ represents a density of the velocity of a representative particle (among infinitely many). The potential $V$ represents the friction and $W$ the inelastic collisions between particles with different velocities. In a different context, equations of this type were also considered in the recent survey \cite{David1} where it is used to describe such phenomena as aggregation, diffusion, and self-propulsion mechanisms among particles.

Several works study the long-time behaviour of solutions to this
equation, including the existence of stationary measures and rates
of convergence to equilibrium \cite{Vil1,Vil2,Mal1}, or develop
numerical schemes for its approximation \cite{Bla1,Mal2}.
The question of existence, however, is not always treated in full
detail in these works. In particular, \cite{Vil1} provides a detailed
existence proof on bounded domains and outlines the argument
for $\R^d$. Under stronger assumptions on $V$ and $W$, existence
in the convex regime can also be established using gradient flow
theory \cite{Amb1,Fig1}, whereas our approach applies to a broader
class of potentials than those covered by these results.

The problem of existence for linear diffusion ($m=1$) was rigorously solved in \cite{Cat1} using a probabilistic technique called propagation of chaos that has been known since \cite{Kac1,Mel1,Szn1} and surveyed in \cite{Cha1,Cha2}. The core of this approach is to consider a system of $N$ particles $X_t^{1,N},\ldots, X_t^{N,N}$ that are governed by $N$ linear SDEs and consider their empirical measure. If it converges to a deterministic limit as $N\to\infty$, we call it propagation of chaos. The name comes from the observation that, under natural assumptions, this property is equivalent to the asymptotic independence of the particles. More precisely, if we fix $k\in\N$, then $\text{Law}(X^{1,N}_t,\ldots,X^{k,N}_t)\xrightarrow{N\to\infty}f_t\otimes\ldots\otimes f_t$ in the weak sense, where $f_t$ is the limit of the empirical measures.

In \cite{Cat1}, the authors considered the following system of SDEs for $i=1,\ldots,N$ and $t>0$:
\begin{equation}\label{linear SDE}
    \begin{cases}
        dX_t^{i,N}=-\left(\nabla V(X_t^{i,N})+\frac{1}{N}\sum_{j=1}^N\nabla W(X_t^{i,N}-X_t^{j,N})\right)dt+\sqrt{2}dB_t^i\\
     X_0^i\sim\mu_0,
    \end{cases}
\end{equation}
where $B_t^1,\ldots, B_t^N$ are independent $d$-dimensional standard Brownian motions. They showed that this problem is well-posed, the sequence of distributions of their empirical measures, i.e,
\[
\text{Law}\left(\mu_\cdot^N\right):=\text{Law}\left(\frac{1}{N}\sum_{k=1}^N\delta_{X_\cdot^{k,N}}\right)\in\mathcal{P}(\mathcal{P}(C([0,T],\R^d))),
\]
is tight and the limit is $\delta_Q\in\mathcal{P}(\mathcal{P}(C([0,T],\R^d)))$, where $Q$ is the law of a solution to \ref{SDE} for $m=1$. For $t>0$, the time marginals of $Q$ have densities $u_t$ solving \ref{granular medium equation}. Moreover $X_\cdot^{1,N}\to Y_\cdot$ in law, where $Y_\cdot$ solves \ref{SDE} for linear diffusion case, i.e. with $m=1$.

Our goal is to adapt this approach to obtain an analogous result
for $m>1$. A key difficulty in constructing a particle approximation
of this SDE is that the diffusion coefficient depends on pointwise
values of the probability density. Since empirical measures are
atomic, they cannot be used directly to evaluate this coefficient.
We address this difficulty by regularizing the empirical measure
and using the resulting density approximation to define the
diffusion coefficients of the particle system.

When $V=W=0$, equation \eqref{granular medium equation} reduces
to the classical porous medium equation
\begin{align}\label{PME}
    \left\{
    \begin{array}{ll}
        \partial_t u=\Delta(u^m)
        & \text{in }\R^d\times(0,T),\\
        u(\cdot,0)=u_0
        & \text{in }\R^d.
    \end{array}
    \right.
    \tag{PME}
\end{align}
Stochastic particle approximations of this equation have been studied
in several settings. In one dimension, Jourdain \cite{Ben1}
established a propagation of chaos result for a formulation in which
the solution represents a cumulative distribution function.
This result was extended to broader classes of functions
in \cite{roc1}. Philipowski \cite{Phi1} constructed an interacting
diffusion approximation for the quadratic porous medium equation, $m=2$,
in arbitrary spatial dimension. Figalli and Philipowski \cite{Fig2}
subsequently extended this approach to general exponents $m>1$,
using an additional approximation of the nonlinearity when $1<m<2$.
Their construction represents nonlinear diffusion through
a regularized interaction drift and includes auxiliary Brownian
noise. In our construction, nonlinear diffusion is instead
approximated through a diffusion coefficient depending on
a regularized empirical density.

Alongside these stochastic approaches, deterministic particle
approximations of the porous medium equation have also been studied.
In \cite{Ama1,Car3}, convergence of such approximations to solutions
of \eqref{PME} was established using a gradient flow approach.
Related results in \cite{Car1,Car2} concern the convergence of
solutions to non-local diffusion equations, whose diffusion
operators involve an interaction kernel, to solutions of
\eqref{granular medium equation}. These results cover the cases
$m\leq1$ and $m=2$, under the assumption that $V$ and $W$ are
semiconvex and grow at most quadratically. In particular, the
authors of \cite{Car3} used the JKO scheme
(see, e.g., \cite{Amb1}) to construct solutions to the following
non-local equation:
\begin{align}\label{car}
    \partial_t u^\epsilon=\frac{m}{m-1}\div(u^\epsilon\nabla\varphi_\epsilon*(\varphi_\epsilon*u^\epsilon)^{m-1}),
    \tag{NLPME}
\end{align}
where $\varphi_\epsilon(x)=\frac{1}{\epsilon^d}\varphi\big(\frac{x}{\epsilon}\big)$ is a standard mollifier, and showed the convergence of $u^\epsilon$ with $\epsilon\to0$ to the solution of \ref{PME}.

Most of the approximations discussed above can be viewed as variants
of the following particle system, which is deterministic for
$\delta=0$ and stochastic for $\delta>0$:
\begin{equation}
\left\{\begin{aligned}
    \text{d}X_t^{i,N} = &  -\frac{m}{m-1}\int\nabla\varphi_\epsilon(y) \left( \frac{1}{N} \sum_{j=1}^N \varphi_\epsilon(X_t^{i,N} - y - X_t^{j,N}) \right)^{m-1} \, \text{d}y\,\text{d}t \\&+ \sqrt{2\delta}  \, \text{d}B_t^i \\
    X_0^{i,N} \sim & \mu_0
\end{aligned}\right.
\end{equation}

In most of the particle approximations discussed above, nonlinear
diffusion is represented through a regularized interaction drift,
while any additional Brownian noise serves as an auxiliary
regularization that vanishes in the limit. Exceptions include
\cite{Ben1}, which is restricted to one spatial dimension, and
\cite{Ben2}, which assumes a Lipschitz diffusion coefficient.

The key feature of our construction is that nonlinear diffusion
is generated entirely by a Brownian term whose coefficient depends
on a regularized empirical density. Our approach applies in arbitrary
spatial dimension and does not require the limiting diffusion
coefficient to be Lipschitz. Crucially, the density-dependent
stochastic term remains present in the limiting SDE: letting
$\delta$ tend to zero removes only the auxiliary viscosity.
More precisely, for $i=1,\ldots,N$, we consider the particle system
\begin{equation}\label{PSDE}
    \left\{
    \begin{aligned}
        \mathrm{d}X_t^{i,N} ={}&
        -\left(
            \nabla V(X_t^{i,N})
            + \frac{1}{N}\sum_{j=1}^N
            \nabla W(X_t^{i,N}-X_t^{j,N})
        \right)\,\mathrm{d}t\\
        &+\sqrt{
            2\int_{\R^d}
            \left(
                \frac{1}{N}\sum_{j=1}^N
                \varphi_\epsilon(X_t^{i,N}-y-X_t^{j,N})
            \right)^{m-1}
            \varphi_\epsilon(y)\,\mathrm{d}y
            +2\delta
        }\,\mathrm{d}B_t^i,\\
        X_0^{i,N}\sim{}&\mu_0.
    \end{aligned}
    \right.
    \tag{PSDE}
\end{equation}

For fixed $\epsilon$ and $\delta$, we establish convergence of particle system,
as $N\to\infty$. The limiting one-particle distributions are the martingale solution to the following non-local SDE:
\begin{equation}
    \left\{\begin{array}{ll}
        dY^{\epsilon,\delta}_t=-(\nabla V(Y_t^{\epsilon,\delta})+\nabla W*\u_t(Y^{\epsilon,\delta}_t))\,dt\\\quad\quad\qquad+\sqrt{2\Big[\varphi_\epsilon*\big(\varphi_\epsilon*\u_t\big)^{m-1}\Big](Y^{\epsilon,\delta}_t)+2\delta}\,dB_t,\\
    \u_t=\text{Law}(Y^{\epsilon,\delta}_t),\\
    Y^{\epsilon,\delta}_0\sim u_0.
    \end{array}\right.
\end{equation}
By It\^o's formula, its time marginals satisfy the associated PDE in the weak sense:
\begin{align}\label{nonlocal}
    \partial_t \u=\div(\u(\nabla V+\nabla W*\u))+\Delta (\u \ \varphi_\epsilon*(\varphi_\epsilon*\u)^{m-1}+\delta \u).
    \tag{NLGME}
\end{align}
We then use entropy estimates to derive Sobolev bounds, following
the approach of \cite{Car3}. These bounds allow us to pass to the
limit as $\epsilon\to0$ and show that both $\u$ and
$\varphi_\epsilon*\u$ converge to a solution $u^\delta$ of the
following equation:

\begin{align}\label{delta}
    \partial_t u^\delta=\div(u^\delta(\nabla V+\nabla W*u^\delta))+\Delta((u^\delta)^m+\delta u^\delta).
\end{align}

Finally, we pass to the limit as $\delta\to0$ to obtain weak solutions
to \eqref{granular medium equation}. We then apply the superposition
principle \cite{Tre1} to construct corresponding martingale solutions
to \eqref{SDE}.

Having outlined our approach, we now give a precise statement
of our main result.

\subsection*{Assumptions and main result}
Throughout this paper, $d\ge1$, $m>1$, and $T<\infty$. We assume the following conditions on $V$, $W$, and the initial law, with a common exponent $q\ge1$, which denotes the growth order of the gradients. Sobolev potentials are always identified with their $C^1$ representatives.
\begin{description}[style=nextline,leftmargin=12mm]
\item[(V)]
$V\in W^{2,\infty}_{\mathrm{loc}}(\R^d;\R)$. For some finite constants
$C_{1,V},C_{2,V},K_V\ge0$ and all $x,y\in\R^d$,
\[
|\nabla V(x)-\nabla V(y)|
\le C_{1,V}(|x-y|\wedge1)(1+|x|^{q}+|y|^{q})
\]
and
\[
(\nabla V(x)-\nabla V(y))\cdot(x-y)
\ge-K_V|x-y|^2-C_{2,V}.
\]
\item[(W)]
$W\in W^{2,\infty}_{\mathrm{loc}}(\R^d;\R)$ is even.
For some $C_{1,W},C_{2,W},K_W\ge0$ and all $x,y\in\R^d$,
\[
|\nabla W(x)-\nabla W(y)|
\le C_{1,W}(|x-y|\wedge1)(1+|x|^{q}+|y|^{q})
\]
and
\[
(\nabla W(x)-\nabla W(y))\cdot(x-y)
\ge-K_W|x-y|^2-C_{2,W}.
\]
\item[($\mu_0$)]
The initial law has a probability density, $\mu_0=u_0(x)\,d x$, with finite Boltzmann entropy
\[
H(u_0):=\int_{\R^d}u_0(x)\log u_0(x)\,d x\in\R,
\qquad \text{with convention}\ 0\log0:=0.
\]
For some $\eta>0$,
\[
\mu_0\in\mathcal P^{\,2q^2+\eta}(\R^d),
\]
where
\[
\mathcal P^p(\R^d)=\left\{\mu\in\mathcal P(\R^d):\int|x|^p\,d\mu(x)<\infty\right\}.
\]
The initial positions are independent with common law $\mu_0$ and are independent of the independent Brownian motions $B^1,\ldots,B^N$.
\end{description}

Assumptions (V) and (W) allow us to work beyond the convex setting
of the gradient flow results discussed above \cite{Amb1,Fig1}.
Moreover, no coercivity condition of power type is imposed on $V$.
The stochastic particle approximation of the non-local equation
\eqref{nonlocal} also accommodates general initial probability
measures: together with (V) and (W), the moment condition in
($\mu_0$) suffices to obtain the particle moment estimates,
identify limit points, and establish convergence of the empirical
measures. At this stage, neither absolute continuity with respect
to Lebesgue measure nor finite initial entropy is required,
so singular initial laws are also admissible. The finite initial
entropy assumption is used in the subsequent density regularity
estimates at the end of Section~\ref{3}.

However, to show convergence with $\epsilon,\delta\to0$ we need to assume additionally the following inequalities holds almost everywhere $x\in\R^d$
\begin{equation}\label{joint}
\begin{aligned}
\Delta V(x)&\le C\bigl[x\cdot\nabla V(x)+c_V(1+|x|^2)\bigr],\\
\Delta W(x)&\le C\bigl[x\cdot\nabla W(x)+c_W(1+|x|^2)\bigr].
\end{aligned}    
\tag{V',W'}
\end{equation}
which give us the following joint condition we will use in section \ref{4}:
\begin{equation}\label{true}
\begin{aligned}
&\Delta V(x)+\Delta V(y)+2\Delta W(x-y)\\
&\quad\le A\Bigl[
x\cdot\nabla V(x)+y\cdot\nabla V(y)
+(x-y)\cdot\nabla W(x-y)\\
&\hspace{30mm}
+c_V(2+|x|^2+|y|^2)+c_W(1+|x-y|^2)
\Bigr]
\end{aligned}
\end{equation}
for some $A>0$.

\begin{remark}
For clarity of presentation, we have chosen to retain the separate
assumptions \eqref{joint}. These are convenient
sufficient conditions for \eqref{true}, and the proofs use
\eqref{joint} only through \eqref{true}. Consequently, all results
stated under \eqref{joint} remain valid when it is replaced by
\eqref{true}, with all other assumptions unchanged.

The condition \eqref{true} admits a strictly larger class
of admissible pairs of potentials. Indeed, in dimension $d=1$,
consider
\[
V(x)=x^4,
\qquad
W(x)=\int_0^x \sin(s^5)\,ds.
\]
These potentials satisfy assumptions $(V)$ and $(W)$,
as well as \eqref{true}, whereas $W$ does not satisfy the
corresponding inequality in \eqref{joint}.
\end{remark}

Now, we are at position to state the main theorem of this paper.

\begin{theorem}\label{main}
Assume $(V),(W),(\mu_0)$ and \eqref{joint}. Then, there
exists a weak solution of \eqref{SDE} in the sense of Definition \ref{Martingale problem} whose marginal densities
form a weak solution of \eqref{granular medium equation} in the sense of Definition \ref{weak} with
initial datum $u_0$.
\end{theorem}

\begin{remark}
We assume that the initial density has finite 
Boltzmann entropy,
\[
    \int_{\R^d} u_0(x)\log u_0(x)\,\mathrm{d}x \in \R,
\]
but do not require finite initial internal energy. In particular,
$u_0$ need not belong to $L^m(\R^d)$, and we allow
\[
    \int_{\R^d} u_0(x)^m\,\mathrm{d}x=+\infty.
\]
For such initial data, estimates controlled by the initial internal
energy, including the standard energy estimates used in gradient
flow constructions, do not yield finite bounds. Our analysis
instead relies on Boltzmann entropy estimates, which remain
available under the stated assumptions.
\end{remark}

The remainder of this paper is organized as follows.
In Section~\ref{2}, we introduce the definitions and background
material used in our proofs, including weak solutions to PDEs
and SDEs, martingale problems, and Sobolev spaces.
In Section~\ref{3}, we establish well-posedness of the stochastic
particle system \eqref{PSDE} and prove its convergence
to \eqref{nonlocal}. We adapt the arguments of \cite{Cat1}
to accommodate the more involved diffusion coefficient.
In Section~\ref{4}, we prove that the solutions constructed
in Section~\ref{3} converge to weak solutions
of \eqref{granular medium equation} and then obtain corresponding
martingale solutions to \eqref{SDE}.

\section{Preliminaries}\label{2}
\subsection*{Weak solutions}

Throughout this section, $T>0$, $d\geq 1$, and $m>1$.
We write $\mathcal P(E)$ for the set of Borel probability measures on a
Polish space $E$. For an open set $\Omega$, the notation
$C_c^l(\Omega;\R^k)$ denotes the space of $l$ times continuously
differentiable functions with compact support in $\Omega$.
The test functions in $C_c^\infty([0,T]\times\R^d;\R)$ are
smooth up to both time endpoints and have support contained in
$[0,T]\times K$ for some compact set $K\subset\R^d$.
They need not vanish at $t=0$ or $t=T$.
We use the notation $f_t(x)=f(t,x)$.

In the definitions below, $u$, $u^\delta$, and $\u$
are nonnegative Borel functions whose spatial sections are probability
densities. The corresponding curves of probability measures are narrowly
continuous: for every $g\in C_b(\R^d;\R)$, the map
\[
    t\longmapsto\int_{\R^d}g(x)u_t(x)\,dx
\]
is continuous, and the same convention applies to $u^\delta$ and
$\u$. The initial density is prescribed and satisfies
$u_0^\delta=u_0^{\varepsilon,\delta}=u_0$.
Whenever a convolution with $\nabla W$ occurs, it is required to be
well defined almost everywhere with respect to the corresponding
solution measure on $[0,T]\times\R^d$.

We begin with the definition of a weak solution to \eqref{granular medium equation}.

\begin{definition}
\label{weak}
We say that $(u_t)_{t\in[0,T]}$ is a weak solution to \eqref{granular medium equation} if, for every
$f\in C_c^\infty([0,T]\times\R^d;\R)$,
\[
\begin{aligned}
 -\int_0^T\int_{\R^d}u\,\partial_t f\,dx\,dt
 &=\int_{\R^d}f_0u_0\,dx
   -\int_{\R^d}f_Tu_T\,dx\\
 &\quad-\int_0^T\int_{\R^d}
    u(\nabla V+\nabla W*u)\cdot\nabla f\,dx\,dt\\
 &\quad+\int_0^T\int_{\R^d}u^m\Delta f\,dx\,dt,
\end{aligned}
\]
and the following integrability condition holds:
\[
 \int_0^T\int_{\R^d}
 \bigl(|\nabla V+\nabla W*u|u+u^m\bigr)\,dx\,dt<\infty.
\]
\end{definition}

We define weak solutions to \eqref{delta} and \eqref{nonlocal} in the same way.
For the nonlocal equation, we use the mollifiers
\[
 \varphi_\varepsilon(x)=\varepsilon^{-d}\varphi(x/\varepsilon),
 \qquad \varepsilon>0,
\]
where $\varphi\in C_c^\infty(\R^d)$ is nonnegative and even, and
$\int_{\R^d}\varphi(x)\,dx=1$.

\begin{definition}
\label{weak delta}
Let $\delta>0$. We say that $(u_t^\delta)_{t\in[0,T]}$ is a weak solution
to (1.4) if, for every
$f\in C_c^\infty([0,T]\times\R^d;\R)$,
\[
\begin{aligned}
 -\int_0^T\int_{\R^d}u^\delta\partial_t f\,dx\,dt
 &=\int_{\R^d}f_0u_0^\delta\,dx
   -\int_{\R^d}f_Tu_T^\delta\,dx\\
 &\quad-\int_0^T\int_{\R^d}
    u^\delta(\nabla V+\nabla W*u^\delta)\cdot\nabla f\,dx\,dt\\
 &\quad+\int_0^T\int_{\R^d}
    \bigl((u^\delta)^m+\delta u^\delta\bigr)\Delta f\,dx\,dt,
\end{aligned}
\]
and
\[
 \int_0^T\int_{\R^d}
 \bigl(|\nabla V+\nabla W*u^\delta|u^\delta
       +(u^\delta)^m\bigr)\,dx\,dt<\infty.
\]
\end{definition}

\begin{definition}
\label{weak nonlocal}
Let $\varepsilon,\delta>0$. We say that
$(u_t^{\varepsilon,\delta})_{t\in[0,T]}$ is a weak solution to (NLGME)
if, for every $f\in C_c^\infty([0,T]\times\R^d;\R)$,
\[
\begin{aligned}
 -\int_0^T\int_{\R^d}
       \u\partial_t f\,dx\,dt&=\int_{\R^d}f_0u_0^{\varepsilon,\delta}\,dx
       -\int_{\R^d}f_Tu_T^{\varepsilon,\delta}\,dx\\
 &\qquad-\int_0^T\int_{\R^d}
       \u
       (\nabla V+\nabla W*\u)
       \cdot\nabla f\,dx\,dt\\
 &\qquad+\int_0^T\int_{\R^d}
       \u
       \left[\varphi_\varepsilon*
         (\varphi_\varepsilon*\u)^{m-1}
         +\delta\right]\Delta f\,dx\,dt,
\end{aligned}
\]
and
\[
 \int_0^T\int_{\R^d}
 |\nabla V+\nabla W*\u|
 \u\,dx\,dt<\infty.
\]
\end{definition}

For fixed $\varepsilon,\delta>0$, no additional $L^m$ assumption is
needed to make the diffusion term in this definition integrable.
Indeed, the unit mass and nonnegativity of $u_t^{\varepsilon,\delta}$ give
\[
 0\leq\varphi_\varepsilon*
       (\varphi_\varepsilon*u_t^{\varepsilon,\delta})^{m-1}
 \leq\|\varphi_\varepsilon\|_{L^\infty(\R^d)}^{m-1}.
\]
Consequently,
\[
\begin{aligned}
 &\int_0^T\int_{\R^d}
 \u
 \left[\varphi_\varepsilon*
   (\varphi_\varepsilon*\u)^{m-1}
   +\delta\right]\,dx\,dt\\
 &\qquad\leq
 T\left(\|\varphi_\varepsilon\|_{L^\infty(\R^d)}^{m-1}
       +\delta\right).
\end{aligned}
\]

We next define solutions to the stochastic equations, starting with strong
solutions to (PSDE). Fix $N\in\mathbb N$ and $\varepsilon,\delta>0$.
Let $B^1,\ldots,B^N$ be independent standard $d$-dimensional Brownian
motions on a given probability space. The initial random variables
$X_0^{1,N},\ldots,X_0^{N,N}$ are independent, have common law with density
$u_0$, and are independent of the Brownian motions.
We use the usual augmentation of the filtration generated by these
initial variables and the Brownian motions.

\begin{definition}
\label{particles}
A family of continuous adapted processes
$(X_t^{i,N})_{t\in[0,T]}$, $i=1,\ldots,N$, is a strong solution to (PSDE)
if, almost surely, for every $i=1,\ldots,N$ and $t\in[0,T]$,
\[
\begin{array}{l}
 X_t^{i,N}
 =X_0^{i,N}
   -\int_0^t\left[
       \nabla V(X_s^{i,N})
       +\frac1N\sum_{j=1}^N\nabla W(X_s^{i,N}-X_s^{j,N})
     \right]ds\\
 \quad+\int_0^t\left[
       2\int_{\R^d}
       \left(\frac1N\sum_{j=1}^N
         \varphi_\varepsilon(X_s^{i,N}-y-X_s^{j,N})
       \right)^{m-1}
       \varphi_\varepsilon(y)\,\text{d}y+2\delta
     \right]^{1/2}\text{d}B_s^i,
\end{array}
\]
where the drift is absolutely integrable on $[0,T]$, almost surely.
The stochastic integral is well defined since its coefficient is bounded
for fixed $\varepsilon,\delta>0$.
\end{definition}

To construct weak solutions to (SDE), we will use the formulation in terms
of a martingale problem. We first introduce the corresponding diffusion
operator.

\begin{definition}
\label{generator}
Let $b:[0,T]\times\R^d\to\R^d$ and
$a:[0,T]\times\R^d\to\R^{d\times d}$ be Borel maps,
with $a(t,x)$ symmetric and nonnegative semidefinite for every $(t,x)$.
The diffusion operator associated with $a$ and $b$ is
\[
 Lf(t,x)
 =\sum_{i=1}^d b_i(t,x)\frac{\partial f}{\partial x_i}(t,x)
  +\frac12\sum_{i,j=1}^d a_{ij}(t,x)
       \frac{\partial^2 f}{\partial x_i\partial x_j}(t,x).
\]
For the stochastic differential equation
\[
 dX_t=b(t,X_t)\,dt+\sigma(t,X_t)\,dB_t,
\]
the covariance matrix is $a=\sigma\sigma^{\mathsf T}$.
The operator $L$ acts only on the spatial variables.
\end{definition}

For $t\in[0,T]$, let
$e_t:C([0,T];\R^d)\to\R^d$ be the evaluation map,
$e_t(\gamma)=\gamma(t)$.
The canonical filtration is $\mathcal F_t=\sigma(e_s:0\leq s\leq t)$.
For $\eta\in\mathcal P(C([0,T];\R^d))$, we write
$\eta_t=(e_t)_\#\eta$ for its time marginal.
For vectors, $|\cdot|$ denotes the Euclidean norm, and for matrices it
denotes the Hilbert--Schmidt norm.

\begin{definition}
\label{Martingale problem}
A probability measure $\eta\in\mathcal P(C([0,T];\R^d))$ solves
the martingale problem associated with $L$ if
\[
 \int_{C([0,T];\R^d)}
 \int_0^T
 \bigl(|b(t,e_t)|+|a(t,e_t)|\bigr)\,dt\,d\eta<\infty
\]
and, for every $f\in C_c^\infty([0,T]\times\R^d;\R)$,
the process
\begin{equation}
\label{eq:prelim-martingale}
 t\longmapsto
 f(t,e_t)-f(0,e_0)
 -\int_0^t\bigl[\partial_s f(s,e_s)+Lf(s,e_s)\bigr]\,ds
\end{equation}
is a martingale under $\eta$ with respect to $(\mathcal F_t)_{t\in[0,T]}$.
A prescribed initial law $\mu_0$ is imposed by requiring $\eta_0=\mu_0$.
\end{definition}

Taking expectations in \eqref{eq:prelim-martingale} yields the weak
formulation of the corresponding Fokker--Planck equation. We formalize it in the following lemma in the case of the \eqref{nonlocal}.
\begin{lemma}
\label{lem:nonlocal-martingale-fokker-planck}
Let $\eta$ be a martingale solution of (1.3), with initial law $\mu_0$,
and assume that $(\eta_t)_{t\in[0,T]}$ is narrowly continuous.
Then $(\eta_t)_{t\in[0,T]}$ is a weak solution of the associated
Fokker--Planck equation. If $\eta_t(dx)=u_t^{\varepsilon,\delta}(x)\,dx$,
then $u^{\varepsilon,\delta}$ is a weak solution of (NLGME),
with the same initial condition.
\end{lemma}

\begin{proof}
Fix $f\in C_c^\infty([0,T]\times\mathbb R^d;\mathbb R)$.
The boundedness of $f$, its derivatives, and $\varphi_\varepsilon$,
together with the drift integrability in the definition of a martingale
solution, justifies Fubini's theorem.
Taking expectations in the martingale identity, applying Fubini's
theorem, and using $\eta_s=(e_s)_\#\eta$, we obtain
\[
\begin{aligned}
 &\int_{\mathbb R^d}f(t,x)\,d\eta_t(x)
  -\int_{\mathbb R^d}f(0,x)\,d\mu_0(x)\\
 &=\int_{C([0,T];\mathbb R^d)}\int_0^t
   \Bigl\{\partial_s f(s,e_s)\\
 &\qquad-\bigl(\nabla V(e_s)+(\nabla W*\eta_s)(e_s)\bigr)
     \cdot\nabla f(s,e_s)\\
 &\qquad+\Bigl(
     \bigl[\varphi_\varepsilon*
       (\varphi_\varepsilon*\eta_s)^{m-1}\bigr](e_s)+\delta
     \Bigr)\Delta f(s,e_s)\Bigr\}\,ds\,d\eta\\[0.5em]
 &=\int_0^t\int_{C([0,T];\mathbb R^d)}
   \Bigl\{\partial_s f(s,e_s)\\
 &\qquad-\bigl(\nabla V(e_s)+(\nabla W*\eta_s)(e_s)\bigr)
     \cdot\nabla f(s,e_s)\\
 &\qquad+\Bigl(
     \bigl[\varphi_\varepsilon*
       (\varphi_\varepsilon*\eta_s)^{m-1}\bigr](e_s)+\delta
     \Bigr)\Delta f(s,e_s)\Bigr\}\,d\eta\,ds\\[0.5em]
 &=\int_0^t\int_{\mathbb R^d}
   \Bigl\{\partial_s f(s,x)\\
 &\qquad-\bigl(\nabla V(x)+(\nabla W*\eta_s)(x)\bigr)
     \cdot\nabla f(s,x)\\
 &\qquad+\Bigl(
     \bigl[\varphi_\varepsilon*
       (\varphi_\varepsilon*\eta_s)^{m-1}\bigr](x)+\delta
     \Bigr)\Delta f(s,x)\Bigr\}\,d\eta_s(x)\,ds.
\end{aligned}
\]
This is the weak formulation of the associated Fokker--Planck equation.
If $\eta_s(dx)=u_s^{\varepsilon,\delta}(x)\,dx$, substituting this
identity yields the weak formulation of (NLGME), with the same initial
condition.
\end{proof}

\begin{definition}
\label{fokker-planck}
A Borel curve $(\nu_t)_{t\in[0,T]}\subset\mathcal P(\R^d)$ is
a weak solution of
\begin{equation}
\label{eq:prelim-fokker-planck}
 \partial_t\nu=L^*\nu
\end{equation}
if
\[
 \int_0^T\int_{\R^d}
 \bigl(|a(t,x)|+|b(t,x)|\bigr)\,d\nu_t(x)\,dt<\infty
\]
and, for every $f\in C_c^\infty([0,T]\times\R^d;\R)$,
\[
\begin{aligned}
 &\int_0^T\int_{\R^d}
   \bigl[\partial_t f(t,x)+Lf(t,x)\bigr]\,d\nu_t(x)\,dt\\
 &\qquad=\int_{\R^d}f(T,x)\,d\nu_T(x)
          -\int_{\R^d}f(0,x)\,d\nu_0(x).
\end{aligned}
\]
Here $L^*$ denotes the formal adjoint in the spatial variables.
\end{definition}

Under the integrability condition in this definition, a probability-valued
solution has a unique narrowly continuous representative; see
[\cite{Tre1}, Remark~2.3]. We use this representative when specifying its values
at individual times.
The converse passage from a Fokker--Planck solution to a solution of the
martingale problem follows from the superposition principle.

\begin{theorem}[Superposition principle; {\normalfont[\cite{Tre1}, Theorem~2.5]}]
\label{thm:prelim-superposition}
Let $a$ and $b$ satisfy the assumptions of
Definition~\ref{generator}, and let
$(\nu_t)_{t\in[0,T]}\subset\mathcal P(\R^d)$ be a narrowly
continuous weak solution of \eqref{eq:prelim-fokker-planck} in the sense
of Definition~\ref{fokker-planck}.
There exists $\eta\in\mathcal P(C([0,T];\R^d))$ solving the
martingale problem associated with $L$ such that
\[
 (e_t)_\#\eta=\nu_t
 \qquad\text{for every }t\in[0,T].
\]
\end{theorem}
For a fixed weak solution $u$ to (PDE), the associated linear operator is
\[
 Lf(t,x)
 =-\bigl(\nabla V(x)+(\nabla W*u_t)(x)\bigr)\cdot\nabla f(t,x)
   +u_t(x)^{m-1}\Delta f(t,x).
\]
Its coefficients are $b=-(\nabla V+\nabla W*u)$ and
$a=2u^{m-1}I_d$, where $I_d$ is the identity matrix.
The integrability required by Theorem~\ref{thm:prelim-superposition}
is exactly the integrability in Definition~\ref{weak}, up to
the constant factor $2\sqrt d$ in the matrix norm.
Thus the theorem gives a martingale solution with time marginals
$u_t(x)\,dx$.
By the representation of continuous local martingales as stochastic
integrals, on a possibly enlarged probability space this yields a weak
solution to (SDE) satisfying
\[
 \operatorname{Law}(Y_t)(dx)=u_t(x)\,dx.
\]
This existence statement does not assert uniqueness of the resulting
path law.

\subsection*{Sobolev spaces}

We recall the Sobolev spaces and embedding results used below.
For the integer-order spaces, we follow \cite{Evans}.

\begin{definition}
\label{sobolev}
Let $k\in\mathbb N$, $1\leq p\leq\infty$, and let
$\Omega\subset\R^d$ be open.
The Sobolev space $W^{k,p}(\Omega)$ consists of the functions
$u\in L^p(\Omega)$ whose weak derivatives $D^\alpha u$ belong to
$L^p(\Omega)$ for every multi-index $\alpha$ with $|\alpha|\leq k$.
Functions that agree almost everywhere are identified.
We equip this space with the norm
\[
 \|u\|_{W^{k,p}(\Omega)}
 =\begin{cases}
 \displaystyle
 \left(\sum_{|\alpha|\leq k}
       \int_\Omega |D^\alpha u(x)|^p\,dx\right)^{1/p},
       &1\leq p<\infty,\\
 \displaystyle
 \sum_{|\alpha|\leq k}
 \operatorname*{ess\,sup}_{x\in\Omega}|D^\alpha u(x)|,
       &p=\infty.
 \end{cases}
\]
\end{definition}

\begin{remark}
For $p=2$, we use the notation $H^k(\Omega)=W^{k,2}(\Omega)$.
\end{remark}

The following Sobolev inequality will be
used in the compactness arguments.

\begin{theorem}[Gagliardo--Nirenberg--Sobolev inequality]
\label{thm:prelim-gns}
Let $d\geq2$ and $1\leq p<d$, and set $p^*=dp/(d-p)$.
There exists a constant $C$, depending only on $d$ and $p$, such that
\[
 \|u\|_{L^{p^*}(\R^d)}
 \leq C\|\nabla u\|_{L^p(\R^d)}
 \qquad\text{for every }u\in W^{1,p}(\R^d).
\]
In dimension $d=1$, the endpoint estimate is
\[
 \|u\|_{L^\infty(\R)}
 \leq\|u'\|_{L^1(\R)}
 \qquad\text{for every }u\in W^{1,1}(\R).
\]
\end{theorem}

We will also use the inhomogeneous Sobolev embeddings.
If $p=d>1$, then for every $p\leq q<\infty$,
\[
 \|u\|_{L^q(\R^d)}
 \leq C(d,q)\|u\|_{W^{1,p}(\R^d)}.
\]
If $d<p<\infty$, then
\[
 \|u\|_{L^\infty(\R^d)}
 \leq C(d,p)\|u\|_{W^{1,p}(\R^d)}.
\]

We also need fractional Sobolev spaces. We use the definition and compact
embedding result from \cite{Nez1}.

\begin{definition}
\label{fractional}
Let $s\in(0,1)$, $1\leq p<\infty$, and let
$\Omega\subset\R^d$ be open. We define
\[
 W^{s,p}(\Omega)
 =\left\{u\in L^p(\Omega):
   \int_\Omega\int_\Omega
   \frac{|u(x)-u(y)|^p}{|x-y|^{d+sp}}\,dx\,dy<\infty
  \right\}.
\]
This is a Banach space with norm
\[
 \|u\|_{W^{s,p}(\Omega)}
 =\left(
   \int_\Omega |u(x)|^p\,dx
   +\int_\Omega\int_\Omega
     \frac{|u(x)-u(y)|^p}{|x-y|^{d+sp}}\,dx\,dy
  \right)^{1/p}.
\]
\end{definition}

\begin{theorem}[Fractional compact embedding; {\normalfont[\cite{Nez1}, Theorem~7.1]}]
\label{thm:prelim-fractional-compactness}
Let $s\in(0,1)$, $1\leq p<\infty$, and $1\leq q\leq p$.
If $\Omega\subset\R^d$ is bounded and open with $C^1$ boundary,
then
\[
 W^{s,p}(\Omega)\subset\subset L^q(\Omega).
\]
\end{theorem}

\subsection*{Weak convergence of probability measures}

We finish by recalling the notions of weak convergence and tightness.
Here $E$ is a Polish space, that is, a separable completely metrizable
space, and $\mathcal P(E)$ denotes its Borel probability measures.
For $a>0$, we write
\[
 \mathcal P_a(\R^d)
 =\left\{\mu\in\mathcal P(\R^d):
     \int_{\R^d}|x|^a\,d\mu(x)<\infty
   \right\}.
\]

\begin{definition}
A sequence $(\mu_n)_{n\in\mathbb N}\subset\mathcal P(E)$ converges
weakly, or narrowly, to $\mu\in\mathcal P(E)$ if
\[
 \lim_{n\to\infty}\int_E f\,d\mu_n=\int_E f\,d\mu
 \qquad\text{for every }f\in C_b(E;\R).
\]
\end{definition}

The compactness criterion for this convergence is expressed in terms of
tightness.

\begin{definition}
A family $(\mu_i)_{i\in I}\subset\mathcal P(E)$ is tight if, for every
$\varepsilon>0$, there exists a compact set $K\subset E$ such that
\[
 \sup_{i\in I}\mu_i(E\setminus K)<\varepsilon.
\]
\end{definition}

\begin{theorem}[Prokhorov]
\label{thm:prelim-prokhorov}
For a Polish space $E$, a family of probability measures in
$\mathcal P(E)$ is tight if and only if it is relatively compact in the
topology of weak convergence.
Equivalently, every sequence of measures chosen from the family has a
subsequence converging weakly to a probability measure on $E$.
\end{theorem}

On $\mathcal P(\R^d)$, weak convergence can be metrized by the
bounded Lipschitz distance.

\begin{definition}
For $\mu,\nu\in\mathcal P(\R^d)$, define
\[
\begin{aligned}
 d_{\mathrm{BL}}(\mu,\nu)
 =\sup\Biggl\{&\left|\int_{\R^d}f\,d(\mu-\nu)\right|:
       f\in C_b(\R^d;\R),\\
 &\|f\|_{L^\infty(\R^d)}\leq1 \,\,\text{and} \sup_{\substack{x,y\in\R^d\\x\ne y}}
       \frac{|f(x)-f(y)|}{|x-y|}\leq1
 \Biggr\}.
\end{aligned}
\]
\end{definition}

\begin{theorem}
\label{thm:prelim-bounded-lipschitz}
The bounded Lipschitz distance is a metric on $\mathcal P(\R^d)$
and induces the topology of weak convergence. In particular,
\[
 d_{\mathrm{BL}}(\mu_n,\mu)\longrightarrow0
 \quad\Longleftrightarrow\quad
 \mu_n\longrightarrow\mu\ \text{weakly}.
\]
\end{theorem}

\section{Stochastic particle approximation}\label{3}
\subsection{Well-posedness of particle systems}

In this section we study the particle systems (PSDE) for fixed $\varepsilon,\delta>0$. We prove that every limit point of the laws of their empirical path measures is concentrated on solutions of the nonlinear martingale problem associated with (NLGME). If this nonlinear martingale problem is unique in the moment class obtained below, this gives propagation of chaos. The argument follows \cite{Cat1}; we present it here to make the paper self-contained. The proof is divided into three steps:
\begin{enumerate}
\item proving that the particle systems (PSDE) are well posed;
\item proving weak compactness of the laws of their empirical measures;
\item identifying their limit points through the nonlinear martingale problem.
\end{enumerate}

All estimates in this section are on the finite interval $[0,T]$; their constants may depend on $T,\varepsilon,\delta$.

To obtain well-posedness, we use the standard Lyapunov criterion, as in \cite[Theorem~3.5]{Kha1}. We recall the version needed here.

\begin{theorem}\label{chaos:thm:Lyapunov}
Let $b:[0,T]\times\R^k\to\R^k$ and $\sigma:[0,T]\times\R^k\to\R^{k\times l}$ be continuous and locally Lipschitz in space, uniformly on compact time intervals. Suppose that a function $f\in C^{1,2}([0,T]\times\R^k)$ satisfies,
\[
\lim_{R\to\infty}\inf_{0\le s\le T,\,|x|\ge R}f(s,x)=\infty,
\]
and, for a finite constant $C$,
\[
\partial_sf(s,x)+b(s,x)\cdot\nabla f(s,x)
+\frac12\operatorname{Tr}\!\left(\sigma(s,x)\sigma(s,x)^{\mathsf T}D^2f(s,x)\right)
\le Cf(s,x).
\]
If $X_0$ is independent of the $l$-dimensional Brownian motion and $\E f(0,X_0)<\infty$, then
\[
X_t=X_0+\int_0^tb(s,X_s)\,d s+\int_0^t\sigma(s,X_s)\,d B_s
\]
has a unique nonexplosive strong solution on $[0,T]$, with continuous paths.
\end{theorem}

The function $f$ is a Lyapunov function for the SDE. To use this theorem, rewrite the particle system in the following vector notation:
\[
X_t^N=\begin{bmatrix}X_t^{1,N}\\ \vdots\\X_t^{N,N}\end{bmatrix},
\qquad
\widetilde V(x)=\begin{bmatrix}\nabla V(x_1)\\ \vdots\\\nabla V(x_N)\end{bmatrix},
\]
\[
\widetilde W(x)=\begin{bmatrix}
N^{-1}\sum_{j=1}^N\nabla W(x_1-x_j)\\
\vdots\\
N^{-1}\sum_{j=1}^N\nabla W(x_N-x_j)
\end{bmatrix},
\]
where $x=(x_1,\ldots,x_N)\in\R^{dN}$. The diffusion matrix is diagonal, with 
\[
\sigma_{ii}(x)=
\left[2\delta+2\int_{\R^d}
\left(\frac1N\sum_{j=1}^N\varphi_\varepsilon(x_i-y-x_j)\right)^{m-1}
\varphi_\varepsilon(y)\,dy\right]^{1/2} I_d.
\]
With the concatenated $dN$-dimensional Brownian motion $B$, the equation becomes
\begin{equation}\label{chaos:eq:vector-system}
\,d X_t^N=-\bigl(\widetilde V(X_t^N)+\widetilde W(X_t^N)\bigr)\,d t
+\sigma(X_t^N)\,d B_t.
\end{equation}

\begin{proposition}\label{chaos:prop:wellposed}
For every $N$ and fixed $\varepsilon,\delta>0$, under (V), (W), and ($\mu_0$), there exists a unique global strong solution of \eqref{chaos:eq:vector-system}.
\end{proposition}
\begin{proof}
The assumptions on the gradients give local Lipschitz continuity of
$\widetilde V$ and $\widetilde W$. The diffusion coefficient is Lipschitz because square root is Lipschitz on the interval $\big[2\delta,2\delta+2||\varphi_\epsilon||_\infty^{m-1}\big]$ and the internal function is Lipschitz and takes values in this interval.

Apply Theorem~\ref{chaos:thm:Lyapunov} with $f(x)=1+\tfrac12|x|^2$. Its generator is
\[
Lf(x)=-x\cdot\widetilde V(x)-x\cdot\widetilde W(x)
+\frac12\operatorname{Tr}\!\left(\sigma(x)\sigma(x)^{\mathsf T}\right).
\]
We estimate the three terms in this order. Taking $y=0$ in (V) gives
\[
x_i\cdot\nabla V(x_i)
\ge-K_V|x_i|^2-C_{2,V}+x_i\cdot\nabla V(0)
\ge-C(1+|x_i|^2).
\]
Consequently,
\[
x\cdot\widetilde V(x)=\sum_{i=1}^Nx_i\cdot\nabla V(x_i)
\ge-C(N+|x|^2).
\]
For the interaction, evenness of $W$ implies oddness of $\nabla W$ and $\nabla W(0)=0$. Symmetrization gives
\begin{align*}
x\cdot\widetilde W(x)
&=\frac1N\sum_{i,j}x_i\cdot\nabla W(x_i-x_j)\\
&=\frac1{2N}\sum_{i,j}(x_i-x_j)\cdot\nabla W(x_i-x_j)\\
&\ge-\frac1{2N}\sum_{i,j}\bigl(C_{2,W}+K_W|x_i-x_j|^2\bigr)
\ge-C(N+|x|^2),
\end{align*}
where $\sum_{i,j}|x_i-x_j|^2\le2N\sum_i|x_i|^2$.
Finally,
\begin{align*}
\frac12\operatorname{Tr}\!\left(\sigma(x)\sigma(x)^{\mathsf T}\right)
&=d\sum_{i=1}^N\left[\delta+
\int\left(\frac1N\sum_{j=1}^N\varphi_\varepsilon(x_i-y-x_j)\right)^{m-1}\varphi_\varepsilon(y)\,d y\right]\\
&\le Nd\bigl(\delta+\|\varphi_\varepsilon\|_\infty^{m-1}\bigr).
\end{align*}
Thus $Lf\le C_{N,\varepsilon,\delta}f$. The function $f$ is coercive and has finite initial expectation. Theorem~\ref{chaos:thm:Lyapunov} gives existence and uniqueness on every finite time interval and hence nonexplosion.
\end{proof}

\subsection{Tightness}

This section proves compactness of the particle approximation. We introduce the empirical path measures
\begin{equation}\label{chaos:eq:empirical}
\mu^N=\frac1N\sum_{k=1}^N\delta_{X_\cdot^{k,N}}
\in\mathcal P(C([0,T];\R^d)).
\end{equation}
We will show that their laws are tight in $\mathcal P(\mathcal P(C([0,T];\R^d)))$. Their limit points will be identified in the next subsection. To establish tightness, we use the exchangeability criterion recalled in \cite[Lemma~4.5]{Mel1}. First, we recall the definition.

\begin{definition}
A system of random particles $X_1,\ldots,X_N$ is exchangeable if, for every permutation $\pi$ of $\{1,\ldots,N\}$,
\[
\operatorname{Law}(X_1,\ldots,X_N)=\operatorname{Law}(X_{\pi(1)},\ldots,X_{\pi(N)}).
\]
\end{definition}

The particle paths in \eqref{PSDE} are exchangeable: the initial variables and Brownian motions have permutation-invariant joint laws, and the coefficients respect the same permutations. Pathwise uniqueness then gives permutation invariance of the solution law. Now we can use crucial result obtained in \cite[Lemma~4.5]{Mel1}.

\begin{lemma}\label{chaos:lem:exchangeable-tightness}
For exchangeable systems, the sequence $\operatorname{Law}(\mu^N)$ is tight in
$\mathcal P(\mathcal P(C([0,T];\R^d)))$ if and only if $\operatorname{Law}(X_\cdot^{1,N})$ is tight in $\mathcal P(C([0,T];\R^d))$.
\end{lemma}

It remains to prove tightness for continuous processes. We use the following form of Kolmogorov's criterion.

\begin{theorem}[Kolmogorov's criterion]
Let $X^n$ be continuous $\R^d$-valued processes on $[0,T]$. Suppose that their initial laws are tight and, for some $p>0$, $\beta>0$, and $C<\infty$,
\[
\sup_n\E|X_t^n-X_s^n|^p\le C|t-s|^{1+\beta},\qquad s,t\in[0,T].
\]
Then their laws are tight in $\mathcal P(C([0,T];\R^d))$.
\end{theorem}

To obtain this increment estimate, we need moment bounds for $X^{1,N}$. We begin with pair differences. In the remainder of this section constants may depend on the indicated initial moments, $T,p,d,V,W,m,\varepsilon,\delta$, but not on $N$.

\begin{lemma}\label{chaos:lem:pair-moments}
Assume (V) and (W), and let $p\ge2$ with $\mu_0\in\mathcal P^p(\R^d)$. The solutions of \eqref{PSDE} satisfy
\[
\sup_{N\ge2}\sup_{0\le t\le T}\E|X_t^{1,N}-X_t^{2,N}|^p\le C.
\]
\end{lemma}
\begin{proof}
Stop the system at $\tau_R=\inf\{t\ge0:|X_t^N|\ge R\}$. If $|X_0^N|\ge R$, take $\tau_R=0$. Apply It\^o's formula to $|X^{1,N}-X^{2,N}|^p$. The stopped stochastic integral has zero expectation: its integrand vanishes after $\tau_R$ and is bounded before $\tau_R$. The formula for the expectation is
\begin{align*}
&\E|X_{t\wedge\tau_R}^{1,N}-X_{t\wedge\tau_R}^{2,N}|^p
=\E|X_0^{1,N}-X_0^{2,N}|^p\\
&\quad-p\E\int_0^{t\wedge\tau_R}|X_r^{1,N}-X_r^{2,N}|^{p-2}
(X_r^{1,N}-X_r^{2,N})
\\
&\hspace{24mm}\cdot
\bigl(\nabla V(X_r^{1,N})-\nabla V(X_r^{2,N})\bigr)\,dr\\
&\quad-\frac pN\sum_{j=1}^N\E\int_0^{t\wedge\tau_R}|X_r^{1,N}-X_r^{2,N}|^{p-2}
(X_r^{1,N}-X_r^{2,N})\cdot\\
&\hspace{45mm}\bigl(\nabla W(X_r^{1,N}-X_r^{j,N})-
\nabla W(X_r^{2,N}-X_r^{j,N})\bigr)\,d r\\
&\quad+p(p+d-2)\E\int_0^{t\wedge\tau_R}|X_r^{1,N}-X_r^{2,N}|^{p-2}
\Biggl[2\delta\\
&\hspace{22mm}+\int\left(\frac1N\sum_{j=1}^N\varphi_\varepsilon(X_r^{1,N}-y-X_r^{j,N})\right)^{m-1}\varphi_\varepsilon(y)\,d y\\
&\hspace{22mm}+\int\left(\frac1N\sum_{j=1}^N\varphi_\varepsilon(X_r^{2,N}-y-X_r^{j,N})\right)^{m-1}\varphi_\varepsilon(y)\,d y
\Biggr]\,d r.
\end{align*}
The diffusion term has this form because $B^1$ and $B^2$ have zero mutual quadratic covariation.

By (V), the first drift integrand, including its negative sign, is at most
\[
pC_{2,V}|X_r^{1,N}-X_r^{2,N}|^{p-2}
+pK_V|X_r^{1,N}-X_r^{2,N}|^p.
\]
For each $j$, apply (W) to $X_r^{1,N}-X_r^{j,N}$ and $X_r^{2,N}-X_r^{j,N}$. Their difference is $X_r^{1,N}-X_r^{2,N}$, so the interaction has the same upper bound with $C_{2,W},K_W$ in place of $C_{2,V},K_V$ after averaging over $j$. Finally, the bracketed diffusion coefficient is bounded by
$2\delta+2\|\varphi_\varepsilon\|_\infty^{m-1}$.

Using $z^{p-2}\le1+z^p$ for $z\ge0$, we obtain
\begin{align*}
\E|X_{t\wedge\tau_R}^{1,N}-X_{t\wedge\tau_R}^{2,N}|^p
&\le\E|X_0^{1,N}-X_0^{2,N}|^p\\
&\quad+C\int_0^t\left(1+
\E|X_{r\wedge\tau_R}^{1,N}-X_{r\wedge\tau_R}^{2,N}|^p\right)\,d r.
\end{align*}
The initial term is at most $2^p\int|x|^p\,d\mu_0(x)$. Gronwall's inequality gives a bound independent of $R$ and $N$. Nonexplosion implies $\tau_R\to\infty$ almost surely for each fixed $N$. Fatou's lemma therefore proves the claim. The argument applies directly to every real $p\ge2$.
\end{proof}

Using this lemma, we obtain moment bounds for the individual particles.

\begin{theorem}\label{chaos:thm:particle-moments}
Under (V), (W), and ($\mu_0$), for every
\[
0<p\le2q+\frac\eta{q},
\]
the solutions of \eqref{PSDE} satisfy
\[
\sup_N\sup_{0\le t\le T}\E|X_t^{1,N}|^p\le C.
\]
\end{theorem}
\begin{proof}
We begin with $p\ge2$. Apply the stopped It\^o formula to $|X^{1,N}|^p$ and take expectations. This gives
\begin{align*}
&\E|X_{t\wedge\tau_R}^{1,N}|^p
=\E|X_0^{1,N}|^p\\
&\quad-p\E\int_0^{t\wedge\tau_R}|X_r^{1,N}|^{p-2}
X_r^{1,N}\cdot\nabla V(X_r^{1,N})\,d r\\
&\quad-\frac pN\sum_{j=1}^N\E\int_0^{t\wedge\tau_R}|X_r^{1,N}|^{p-2}
X_r^{1,N}\cdot\nabla W(X_r^{1,N}-X_r^{j,N})\,d r\\
&\quad+p(p+d-2)\E\int_0^{t\wedge\tau_R}|X_r^{1,N}|^{p-2}
\Biggl[\delta\\
&\hspace{27mm}+\int\left(\frac1N\sum_{j=1}^N
\varphi_\varepsilon(X_r^{1,N}-y-X_r^{j,N})\right)^{m-1}\varphi_\varepsilon(y)\,d y\Biggr]\,d r.
\end{align*}
The stochastic integral has zero expectation for the same reason as in Lemma~\ref{chaos:lem:pair-moments}. The lower bound for $V$ gives
\[
-p|X_r^{1,N}|^{p-2}X_r^{1,N}\cdot\nabla V(X_r^{1,N})
\le C\bigl(1+|X_r^{1,N}|^p\bigr).
\]
For the interaction, the growth bound and H\"older's inequality give
\begin{align*}
&\E\left[\mathbf{1}_{\{r<\tau_R\}}|X_r^{1,N}|^{p-1}
|\nabla W(X_r^{1,N}-X_r^{j,N})|\right]\\
&\quad\le C\left(\E|X_{r\wedge\tau_R}^{1,N}|^p\right)^{(p-1)/p}
\left[1+\left(\E|X_r^{1,N}-X_r^{j,N}|^{pq}\right)^{1/p}\right].
\end{align*}
The second bracket is bounded uniformly in $r,N,j$ by Lemma~\ref{chaos:lem:pair-moments}, since
\[
pq\le2q^2+\eta.
\]
For $j=1$ the interaction is zero. For the other indices, exchangeability gives the same bound. Young's inequality therefore bounds the averaged interaction contribution by
$C(1+\E|X_{r\wedge\tau_R}^{1,N}|^p)$.
The diffusion term has the same bound, because its coefficient is at most $\delta+\|\varphi_\varepsilon\|_\infty^{m-1}$.

Combining the estimates, and writing
$m_p(t)=\E|X_{t\wedge\tau_R}^{1,N}|^p$ in this proof, we obtain
\[
m_p(t)\le m_p(0)+C\int_0^t(1+m_p(r))\,d r.
\]
Gronwall's inequality gives $m_p(t)\le(m_p(0)+CT)e^{CT}$ for $t\le T$, uniformly in $R,N$. Let $R\to\infty$ and apply Fatou's lemma. This proves the assertion for $p\ge2$. For $0<p<2$, use Jensen's inequality and the bound for $p=2$. The constant in this estimate may grow with $T$.
\end{proof}

\begin{remark}\label{chaos:rem:radial-W}
If $W$ is radially symmetric, the moment estimate can be obtained directly from the corresponding initial moment. In particular, for tightness and identification of limit points, the moment condition in ($\mu_0$) can then be weakened to
\[
\mu_0\in\mathcal P^{\,2q+\eta/q}(\R^d).
\]
Here is the justification under the present assumptions. Write $W(x)=w(|x|)$. Local Lipschitz continuity and $\nabla W(0)=0$ give $w'(r)/r\ge-C$ for $0<r\le1$. For $r\ge1$, (W) with $y=0$ gives the same bound. Hence
$\nabla W(x)=\frac{w'(|x|)}{|x|}x$ has radial coefficient bounded below by $-C$ on $x\ne0$.

For $p\ge2$, monotonicity of $x\mapsto|x|^{p-2}x$ gives
\[
0\le\bigl(|x|^{p-2}x-|y|^{p-2}y\bigr)\cdot(x-y)
\le C_p(|x|^p+|y|^p).
\]
Symmetrizing the interaction contribution to the generator of $N^{-1}\sum_i|x_i|^p$ therefore gives
\begin{align*}
&-\frac p{2N^2}\sum_{i,j}
\bigl(|x_i|^{p-2}x_i-|x_j|^{p-2}x_j\bigr)\cdot\nabla W(x_i-x_j)\\
&\hspace{30mm}\le\frac C{N^2}\sum_{i,j}(|x_i|^p+|x_j|^p)
\le\frac C N\sum_i|x_i|^p.
\end{align*}
Together with the external-potential and diffusion bounds, stopped It\^o's formula and Gronwall's inequality now propagate every finite initial moment of order $p\ge2$. Exchangeability identifies the expected average with $\E|X_t^{1,N}|^p$. The weaker initial condition displayed above thus supplies all the moments used below. The same symmetrization principle is used in the radial case in \cite{Bla1}.
\end{remark}

We are now ready to estimate increments and prove tightness of the one-particle laws.

\begin{theorem}\label{chaos:thm:increments}
Under (V), (W), and ($\mu_0$), for every
\[
2<p\le 2+\frac{\eta}{q^2},
\]
the solutions of \eqref{PSDE} satisfy
\[
\sup_N\E|X_t^{1,N}-X_s^{1,N}|^p
\le C|t-s|^\alpha,\qquad s,t\in[0,T],\alpha\in[1,p/2).
\]
\end{theorem}
\begin{proof}
Fix $0\le s<t\le T$, and use It\^o's formula for
$|X_t^{1,N}-X_s^{1,N}|^p$, where $p>2$, to get
\begin{align*}
&\E|X_t^{1,N}-X_s^{1,N}|^p\\
&\quad\le p\E\int_s^t|\nabla V(X_u^{1,N})|
|X_u^{1,N}-X_s^{1,N}|^{p-1}\,du\\
&\qquad+\frac pN\sum_{j=1}^N\E\int_s^t
|\nabla W(X_u^{1,N}-X_u^{j,N})|
|X_u^{1,N}-X_s^{1,N}|^{p-1}\,du\\
&\qquad+C\int_s^t\E|X_u^{1,N}-X_s^{1,N}|^{p-2}\,du.
\end{align*}
Here we have used boundedness of $\varphi_\epsilon$ to estimate the
second-order term. 

By H\"older's inequality and the growth assumptions, we get
\begin{align*}
&\E|X_t^{1,N}-X_s^{1,N}|^p\\
&\quad\le C\int_s^t
\left[1+\left(\E|X_u^{1,N}|^{pq}\right)^{1/p}\right]
\left(\E|X_u^{1,N}-X_s^{1,N}|^p\right)^{(p-1)/p}\,du\\
&\qquad+C\int_s^t
\left[1+\left(\E|X_u^{1,N}-X_u^{2,N}|^{pq}\right)^{1/p}\right]
\left(\E|X_u^{1,N}-X_s^{1,N}|^p\right)^{(p-1)/p}\,du\\
&\qquad+C\int_s^t
\left(\E|X_u^{1,N}-X_s^{1,N}|^p\right)^{(p-2)/p}\,du.
\end{align*}
In the interaction term we have used exchangeability and
$\nabla W(0)=0$. Since
\[
pq\le2q+\frac\eta q\le2q^2+\eta,
\]
the moments on the right-hand side are bounded uniformly in $N$ and
$u\in[0,T]$ by Lemma~\ref{chaos:lem:pair-moments} and
Theorem~\ref{chaos:thm:particle-moments}.

We write $m_p(t):=\E|X_t^{1,N}-X_s^{1,N}|^p$ and obtain, for
some constant $C$ independent of $N$ and $s$,
\[
m_p(t)\le C\int_s^t
\left(m_p(u)^{(p-1)/p}+m_p(u)^{(p-2)/p}\right)\,du.
\]
By iterating this inequality we get
\begin{align*}
m_p(t)
&\le C\int_s^t
\left(\int_s^u
\left[m_p(v)^{(p-1)/p}+m_p(v)^{(p-2)/p}\right]\,dv
\right)^{(p-1)/p}\,du\\
&\quad+C\int_s^t
\left(\int_s^u
\left[m_p(v)^{(p-1)/p}+m_p(v)^{(p-2)/p}\right]\,dv
\right)^{(p-2)/p}\,du.
\end{align*}
Moreover, by Theorem~\ref{chaos:thm:particle-moments},
\[
m_p(t)\le2^{p-1}
\left(\E|X_t^{1,N}|^p+\E|X_s^{1,N}|^p\right)\le C.
\]
Using this fact in the iterated inequality, we obtain
\begin{align*}
m_p(t)
&\le C(t-s)^{1+(p-1)/p}+C(t-s)^{1+(p-2)/p}\\
&\le C_T(t-s)^{1+(p-2)/p},
\end{align*}
since $t-s\le T$ and $p>2$.

We can continue this iteration and obtain, with constant $C$ bigger with each iteration $k$, 
\begin{align*}
    m_p(t)\leq C|t-s|^{1+\frac{p-2}{p}+\left(\frac{p-2}{p}\right)^2+\ldots+\left(\frac{p-2}{p}\right)^k}
\end{align*}
\end{proof}

\begin{corollary}\label{chaos:cor:tightness}
The sequence $\operatorname{Law}(\mu^N)$ is tight in
$\mathcal P(\mathcal P(C([0,T];\R^d)))$.
\end{corollary}
\begin{proof}
Choose $p>2$ in the range of Theorem~\ref{chaos:thm:increments}.
Then $p/2>1$. The initial laws all equal $\mu_0$ and are tight.
Kolmogorov's criterion gives tightness of
$\operatorname{Law}(X_\cdot^{1,N})$, and
Lemma~\ref{chaos:lem:exchangeable-tightness} gives the assertion.
\end{proof}

\subsection{Martingale problem}

In this section we show that every limit point of $\operatorname{Law}(\mu^N)$ is concentrated on solutions of the nonlinear martingale problem associated with (NLGME). Denote such a limit by
\[
\Pi_\infty^{\varepsilon,\delta}
\in\mathcal P(\mathcal P(C([0,T];\R^d))).
\]
Thus $\Pi_\infty^{\varepsilon,\delta}$ is a probability measure on path laws. We will prove that for $\Pi_\infty^{\varepsilon,\delta}$-almost every $P$, the initial marginal is $P_0=\mu_0$ and, for every $\psi\in C_c^2(\R^d)$,
\begin{equation}\label{chaos:eq:nonlinear-martingale}
\begin{aligned}
M_t^{\psi,P}
&=\psi(Y_t)-\psi(Y_0)
+\int_0^t\nabla\psi(Y_r)\cdot
\bigl(\nabla V(Y_r)+(\nabla W*P_r)(Y_r)\bigr)\,d r\\
&\quad-\int_0^t\Delta\psi(Y_r)
\bigl[\delta+\varphi_\varepsilon*(\varphi_\varepsilon*P_r)^{m-1}(Y_r)\bigr]\,d r
\end{aligned}
\end{equation}
is a $P$-martingale. Here $Y_t=e_t$ is the canonical process, $P_t$ is its time marginal, and the filtration is
$\mathcal F_t=\sigma(Y_r:0\le r\le t)$, with the usual completion when needed. The superscript $P$ in \eqref{chaos:eq:nonlinear-martingale} records that the coefficients depend on the marginal curve of the law under consideration.

To prove this, fix $k\in\mathbb N$, $\psi\in C_c^2(\R^d)$, $g\in C_b((\R^d)^k)$, and $0\le s_1\le\cdots\le s_k\le s<t\le T$. For path laws $Q$ satisfying
\[
\int_0^T\!\int|x|^{q}\,d Q_r(x)\,d r<\infty,
\]
define the functional
\[
F(Q)=\left\langle Q,
\bigl(M_t^{\psi,Q}-M_s^{\psi,Q}\bigr)
g(Y_{s_1},\ldots,Y_{s_k})\right\rangle.
\]
This moment assumption ensures integrability of the interaction term. All the other terms are bounded after multiplication by the spatial derivatives of $\psi$. It will be enough to prove
\[
\int F(Q)^2\,d\Pi_\infty^{\varepsilon,\delta}(Q)=0
\]
for every such martingale test, and then use a countable determining family. The required moment bounds under the limiting law will be verified below.

Before passing to the limit, we establish the appropriate continuity statement. We must retain the factor $g(Y_{s_1},\ldots,Y_{s_k})$ throughout: in general, the functional involves the path law and cannot be rewritten using only its one-time marginals.

\begin{lemma}\label{chaos:lem:functional-continuity}
Suppose that $Q_n\rightharpoonup Q$ in $\mathcal P(C([0,T];\R^d))$ and, for some $\kappa>0$,
\[
\sup_n\int_0^T\!\int|x|^{q+\kappa}\,d Q_{n,r}(x)\,d r<\infty.
\]
Then $F(Q_n)\to F(Q)$.
\end{lemma}
\begin{proof}
The endpoint terms are bounded continuous functions of the path. Since $\nabla\psi$ has compact support, $\nabla\psi\cdot\nabla V$ is bounded and continuous. Its time integral, multiplied by $g(Y_{s_1},\ldots,Y_{s_k})$, is also a bounded continuous path functional by dominated convergence. The term $\delta\Delta\psi$ is treated in the same way. It remains to consider the two convolutions.

We start with the interaction. The growth assumption gives
\[
\int|\nabla W(x-z)|\,d Q_r(z)
\le C\left(1+|x|^{q}+\int|z|^{q}\,d Q_r(z)\right)
\]
for almost every $r$, and similarly for $Q_n$. The limiting law has the higher time-integrated moment in the statement: the function
\[
\omega\longmapsto\int_0^T|\omega(r)|^{q+\kappa}\,dr
\]
is nonnegative and continuous on path space, so the Portmanteau theorem gives the required lower semicontinuity of its integral.

Let $\chi_R\in C_c^\infty(\R^d)$ be even, between zero and one, equal to one on $B_R$, and zero outside $B_{2R}$. 
For a fixed $R$, write the truncated interaction as an integral against $Q_n\otimes Q_n$. With both path variables in $C([0,T];\R^d)$, it is
\begin{align*}
\int_s^t\!\iint
&g(\omega(s_1),\ldots,\omega(s_k))\,
\nabla\psi(\omega(r))\cdot\nabla W(\omega(r)-\widetilde\omega(r))\\
&\hspace{20mm}\times\chi_R(\omega(r)-\widetilde\omega(r))
\,d Q_n(\omega)\,d Q_n(\widetilde\omega)\,d r.
\end{align*}
For each $r$, the integrand is a bounded continuous function of the two paths, with a bound independent of $r$ and $n$. Since $Q_n\otimes Q_n\rightharpoonup Q\otimes Q$, it converges after integration in the path variables. Dominated convergence then passes the time integral to the limit. This proves continuity of the truncated interaction and keeps the full past-dependent factor $g$.

We now control the removed tail. For $R\ge1$, the growth assumption gives
\[
|\nabla W(x-z)|\mathbf{1}_{\{|x-z|>R\}}
\le CR^{-\kappa}(1+|x|^{q+\kappa}+|z|^{q+\kappa}).
\]
Taking absolute values, bounding $g$ and $\nabla\psi$, and using Fubini, we obtain
\[
|F(Q_n)-F_R(Q_n)|
\le CR^{-\kappa}\left(T+\int_0^T\!\int|x|^{q+\kappa}\,d Q_{n,r}(x)\,d r\right).
\]
The same estimate holds for $Q$. Thus the interaction tails vanish uniformly in $n$ as $R\to\infty$. Convergence for fixed $R$, followed by $R\to\infty$, proves convergence of the full interaction. The uniform tail estimate justifies integration against the changing measures $Q_n$.

We next treat the nonlinear diffusion. Set
\begin{align*}
H_{n,r}(x)&=\Delta\psi(x)\,[\varphi_\varepsilon*(\varphi_\varepsilon*Q_{n,r})^{m-1}](x),\\
H_r(x)&=\Delta\psi(x)\,[\varphi_\varepsilon*(\varphi_\varepsilon*Q_r)^{m-1}](x).
\end{align*}
For each fixed $r$, weak convergence on path space implies $Q_{n,r}\rightharpoonup Q_r$. The functions $\varphi_\varepsilon*Q_{n,r}$ are uniformly bounded by $\|\varphi_\varepsilon\|_\infty$ and uniformly Lipschitz, with Lipschitz constant at most $\|\nabla\varphi_\varepsilon\|_\infty$. They converge pointwise to $\varphi_\varepsilon*Q_r$. Arzel\`a--Ascoli therefore gives uniform convergence on every compact spatial set. The function $z\mapsto z^{m-1}$ is uniformly continuous on $[0,\|\varphi_\varepsilon\|_\infty]$, for every $m>1$. Taking this power preserves local uniform convergence. Because $\varphi_\varepsilon$ has compact support, the outer convolution preserves it as well.

Choose a compact set containing $\operatorname{supp}\psi$. Outside this set $H_{n,r}$ and $H_r$ vanish. Consequently,
\[
\|H_{n,r}-H_r\|_\infty\longrightarrow0
\quad\text{for every fixed }r,
\]
and
\[
\|H_{n,r}\|_\infty+\|H_r\|_\infty
\le2\|\Delta\psi\|_\infty\|\varphi_\varepsilon\|_\infty^{m-1}.
\]
The difference of the diffusion functionals is bounded by
\begin{align*}
&\|g\|_\infty\int_s^t\|H_{n,r}-H_r\|_\infty\,\text{d} r\\
&\quad+\int_s^t
\left|\int H_r(Y_r)g(Y_{s_1},\ldots,Y_{s_k})\,\text{d} Q_n
-\int H_r(Y_r)g(Y_{s_1},\ldots,Y_{s_k})\,\text{d} Q\right|\,\text{d} r.
\end{align*}
The first term tends to zero by dominated convergence. For each fixed $r$, the integrand in the second term tends to zero by weak convergence of the path laws; it is the integral of a bounded continuous function of the path. The displayed uniform bound permits dominated convergence in time. This proves convergence of the diffusion term.

Combining these observations proves that $F_R$ is bounded and continuous on the whole space of path laws. Finally,
\[
|F(Q_n)-F(Q)|\le |F(Q_n)-F_R(Q_n)|
+|F_R(Q_n)-F_R(Q)|+|F_R(Q)-F(Q)|.
\]
Let first $n\to\infty$ and then $R\to\infty$ to obtain the assertion.

The extra control of tails is substantive. For instance, in dimension one take
\[
W(x)=\frac{|x|^{q+1}}{q+1},\qquad
\nu_n=(1-n^{-q})\delta_0+n^{-q}\delta_n.
\]
Then $\nu_n\rightharpoonup\delta_0$ and $\int|x|^{q}\,d\nu_n=1$, but
\[
(\nabla W*\nu_n)(0)=-1,
\qquad (\nabla W*\delta_0)(0)=0.
\]
Thus boundedness of the moment at the growth exponent alone does not ensure convergence. Taking path laws supported on these constant paths, $g=1$, and $\psi'(0)=1$ also gives a failure of continuity of the interaction functional. A higher moment, as in the lemma, provides the required uniform integrability.
\end{proof}

We are now ready to prove the main theorem of this section.

\begin{theorem}\label{chaos:thm:limit-identification}
Every weak limit $\Pi_\infty^{\varepsilon,\delta}$ of $\operatorname{Law}(\mu^N)$ is concentrated on solutions of the nonlinear martingale problem \eqref{chaos:eq:nonlinear-martingale} with initial marginal $\mu_0$.
\end{theorem}
\begin{proof}
We first record the moments available under the limiting law. Fix
\[
0<\kappa<\frac\eta{q}.
\]
Theorem~\ref{chaos:thm:particle-moments} and exchangeability give
\[
\sup_N\E\int_0^T\!\int|x|^{2q+\kappa}\,d\mu_r^N(x)\,d r<\infty.
\]
The corresponding functional of a path law is nonnegative and lower semicontinuous. Hence the Portmanteau theorem gives
\[
\int\left(\int_0^T\!\int|x|^{2q+\kappa}\,d Q_r(x)\,d r\right)
\,d\Pi_\infty^{\varepsilon,\delta}(Q)<\infty.
\]
In particular, the inner integral is finite for almost every $Q$, and all martingale tests are integrable for those laws. The initial empirical measure converges in probability to $\mu_0$ by the law of large numbers. Since $Q\mapsto Q_0$ is continuous, $Q_0=\mu_0$ for $\Pi_\infty^{\varepsilon,\delta}$-almost every $Q$.

It remains to show that $\int F(Q)^2\,d\Pi_\infty^{\varepsilon,\delta}(Q)=0$. We cannot pass an unbounded $F^2$ directly through weak convergence of $\operatorname{Law}(\mu^N)$. We will use the bounded continuous functions $F_R^2$ and a separate squared tail estimate.

For the empirical measure, introduce the particle martingales
\begin{align*}
&M_t^{\psi,i,N}=\psi(X_t^{i,N})-\psi(X_0^{i,N})\\
&\quad+\int_0^t\nabla\psi(X_r^{i,N})\cdot\nabla V(X_r^{i,N})\,dr\\
&\quad+\frac1N\sum_{j=1}^N\int_0^t\nabla\psi(X_r^{i,N})\cdot
\nabla W(X_r^{i,N}-X_r^{j,N})\,dr-\delta\int_0^t\Delta\psi(X_r^{i,N})\,dr\\
&\quad-\int_0^t\Delta\psi(X_r^{i,N})\int_{\R^d}
\left(\frac1N\sum_{j=1}^N\varphi_\varepsilon(X_r^{i,N}-y-X_r^{j,N})\right)^{m-1}
\varphi_\varepsilon(y)\,dy\,dr.
\end{align*}
It\^o's formula identifies them as
\[
M_t^{\psi,i,N}=\int_0^t H_r^{i,N}\, \text{d}B_r^i,
\]
where
\begin{align*}
H_r^{i,N}
&=\nabla\psi(X_r^{i,N})\Biggl[2\delta\\
&\qquad+2\int_{\R^d}
\left(\frac1N\sum_{j=1}^N\varphi_\varepsilon(X_r^{i,N}-y-X_r^{j,N})\right)^{m-1}
\varphi_\varepsilon(y)\,dy\Biggr]^{1/2}.
\end{align*}
Their integrands satisfy
\[
|H_r^{i,N}|^2\le
2\|\nabla\psi\|_\infty^2\bigl(\delta+\|\varphi_\varepsilon\|_\infty^{m-1}\bigr).
\]
Thus these are square-integrable martingales. For $i\ne j$ their mutual\linebreak quadratic covariation vanishes, since the driving Brownian motions are independent. In particular, the increments over $[s,t]$ are conditionally orthogonal given $\mathcal F_s$. The factors $g(X_{s_1}^{i,N},\ldots,X_{s_k}^{i,N})$ are bounded and $\mathcal F_s$-measurable, so the same orthogonality holds after multiplication by these factors.

At the empirical measure,
\[
F(\mu^N)=\frac1N\sum_{i=1}^N
\bigl(M_t^{\psi,i,N}-M_s^{\psi,i,N}\bigr)
g(X_{s_1}^{i,N},\ldots,X_{s_k}^{i,N}).
\]
By exchangeability, expanding the square gives
\begin{align*}
\E F(\mu^N)^2
&=\frac1N\E\left[
\bigl(M_t^{\psi,1,N}-M_s^{\psi,1,N}\bigr)^2
g(X_{s_1}^{1,N},\ldots,X_{s_k}^{1,N})^2\right]\\
&\quad+\frac{N-1}{N}\E\Bigl[
\bigl(M_t^{\psi,1,N}-M_s^{\psi,1,N}\bigr)
\bigl(M_t^{\psi,2,N}-M_s^{\psi,2,N}\bigr)\\
&\hspace{35mm}\times
g(X_{s_1}^{1,N},\ldots,X_{s_k}^{1,N})
g(X_{s_1}^{2,N},\ldots,X_{s_k}^{2,N})\Bigr].
\end{align*}
The second expectation is zero by the conditional orthogonality just proved. The first is bounded, so
\[
\E F(\mu^N)^2\le\frac{C_{\varepsilon,\delta,\psi,g,T}}N.
\]

We next estimate the squared truncation error. Jensen's inequality in the two path variables and in time, followed by the polynomial growth bound, gives, for every admissible path law $Q$,
\begin{align*}
|F(Q)-F_R(Q)|^2
&\le C\int_s^t\!\iint
|\nabla W(x-z)|^2\mathbf{1}_{\{|x-z|>R\}}\,d Q_r(x)\,d Q_r(z)\,d r\\
&\le CR^{-\kappa}\left(T+
\int_0^T\!\int|x|^{2q+\kappa}\,d Q_r(x)\,d r\right).
\end{align*}
For the last step, use
\[
(1+|x-z|^{2q})\mathbf{1}_{\{|x-z|>R\}}
\le CR^{-\kappa}(1+|x|^{2q+\kappa}+|z|^{2q+\kappa}),
\quad R\ge1.
\]
The moment bounds at the start of the proof imply
\begin{align*}
\sup_N\E|F(\mu^N)-F_R(\mu^N)|^2&\le CR^{-\kappa},\\
\int|F(Q)-F_R(Q)|^2\,d\Pi_\infty^{\varepsilon,\delta}(Q)
&\le CR^{-\kappa}.
\end{align*}
Lemma~\ref{chaos:lem:functional-continuity} says that $F_R$ is bounded and continuous on the whole space of path laws. Along the subsequence defining $\Pi_\infty^{\varepsilon,\delta}$, we may therefore pass $F_R^2$ to the limit:
\begin{align*}
\int F_R(Q)^2\,d\Pi_\infty^{\varepsilon,\delta}(Q)
&=\lim_{N\to\infty}\E F_R(\mu^N)^2\\
&\le\limsup_{N\to\infty}
\left(2\E F(\mu^N)^2+2\E|F_R(\mu^N)-F(\mu^N)|^2\right)\\
&\le CR^{-\kappa}.
\end{align*}
It follows that
\begin{align*}
\int F(Q)^2\,d\Pi_\infty^{\varepsilon,\delta}(Q)
&\le2\int|F(Q)-F_R(Q)|^2\,d\Pi_\infty^{\varepsilon,\delta}(Q)\\
&\quad+2\int F_R(Q)^2\,d\Pi_\infty^{\varepsilon,\delta}(Q)
\le CR^{-\kappa}.
\end{align*}
Letting $R\to\infty$ proves the desired identity.

We finish the countability argument. Choose a countable family of spatial tests dense in $C_c^2$ on each compact set, rational time tuples (also including $T$), and countable uniformly dense families of compactly supported continuous cylinder functions $g$, together with the constant function. Intersecting the corresponding full-measure sets gives a single set of path laws on which all these martingale identities hold. For each such law, polynomial integrability implies
\[
\E_Q\int_0^T
\bigl|\nabla\psi(Y_r)\cdot(\nabla V(Y_r)+(\nabla W*Q_r)(Y_r))\bigr|\,d r<\infty.
\]
The diffusion integrand is bounded. Consequently $M_t^{\psi,Q}$ is continuous in $t$ both almost surely and in $L^1(Q)$. Approximation extends the identities from the countable spatial tests and rational times to all tests and times. A monotone-class argument extends the cylinder functions to all bounded $\mathcal F_s$-measurable functions. Thus $M^{\psi,Q}$ is a martingale in the canonical filtration for almost every $Q$.

The global integrability required in Definition~\ref{Martingale problem}
follows from the same moment bound and the growth of the gradients.
For a time-dependent test $f(t,x)=h(t)\psi(x)$, the martingale identity
follows by integration by parts in time. Approximation by finite sums
of such tests, in the $C^{1,2}$ norm on a common compact support,
gives the identity for every
$f\in C_c^\infty([0,T]\times\R^d)$.
This proves the assertion in the precise sense of that definition.
\end{proof}

If we know that the nonlinear martingale problem
has a unique solution $Q^{\varepsilon,\delta}$ among path laws with
finite time-integrated moment of order $2q+\kappa$, then
Theorem~\ref{chaos:thm:limit-identification} implies that every limit
point is $\delta_{Q^{\varepsilon,\delta}}$. Tightness therefore gives
\[
\operatorname{Law}(\mu^N)\rightharpoonup\delta_{Q^{\varepsilon,\delta}}.
\]
In particular, the empirical measures converge in probability to this deterministic path law. Exchangeability then gives, for each fixed $k$,
\[
\operatorname{Law}(X_\cdot^{1,N},\ldots,X_\cdot^{k,N})
\rightharpoonup (Q^{\varepsilon,\delta})^{\otimes k}.
\]
This is exactly propagation of chaos. 

Without uniqueness, Theorem~\ref{chaos:thm:limit-identification} still gives
existence of a nonlinear martingale solution: choose one path law
from the full-measure set in that theorem.

\begin{corollary}\label{chaos:cor:nonlocal-existence}
Assume $(V),(W)$ and $(\mu_0)$. Then for every
$\varepsilon,\delta>0$, equation \eqref{nonlocal} has at least one weak
solution with initial density $u_0$. This solution has narrowly
continuous marginals and, for fixed $\varepsilon,\delta$,
\[
\sqrt{u^{\varepsilon,\delta}}\in
L^2(0,T;H^1(\R^d)).
\]
\end{corollary}
\begin{proof}
Choose a path law $P$ as in Theorem~\ref{chaos:thm:limit-identification},
with finite time-integrated moment of order $2q+\kappa$.
Its marginal curve is narrowly continuous, since its paths are
continuous. Moreover, Jensen's inequality and the growth bounds give
\begin{align*}
&\int_0^T\int_{\R^d}
|\nabla V(x)+(\nabla W*P_r)(x)|^2\,dP_r(x)\,dr\\
&\quad\le C\int_0^T\int_{\R^d}
(1+|x|^{2q})\,dP_r(x)\,dr<\infty.
\end{align*}
The covariance coefficient is bounded. Thus the global integrability
condition in Definition~\ref{Martingale problem} is satisfied.
Lemma~\ref{lem:nonlocal-martingale-fokker-planck} gives the corresponding
probability-measure solution. Explicitly, for every
$\psi\in C_c^2(\R^d)$ and $t\in[0,T]$,
\begin{align*}
&\int_{\R^d}\psi\,dP_t-\int_{\R^d}\psi\,d\mu_0\\
&\quad=-\int_0^t\int_{\R^d}
\nabla\psi(x)\cdot
\bigl(\nabla V(x)+(\nabla W*P_r)(x)\bigr)\,dP_r(x)\,dr\\
&\qquad+\int_0^t\int_{\R^d}\Delta\psi(x)
\Bigl[\delta+\bigl(\varphi_\varepsilon*
(\varphi_\varepsilon*P_r)^{m-1}\bigr)(x)\Bigr]\,dP_r(x)\,dr.
\end{align*}

We check the hypotheses of the density regularity theorem
\cite[Theorem~7.4.1]{Kri1}. The scalar coefficient of $\Delta$ satisfies
\[
\delta\le \delta+\varphi_\varepsilon*
(\varphi_\varepsilon*P_r)^{m-1}
\le\delta+\|\varphi_\varepsilon\|_\infty^{m-1}.
\]
Its spatial gradient is bounded, uniformly in $r$, by
$\|\nabla\varphi_\varepsilon\|_{L^1}
\|\varphi_\varepsilon\|_\infty^{m-1}$; here the derivative is taken
on the outer convolution. The square integrability of the drift was
proved above. The time-integrated second moment is finite, so the
required squared logarithmic moment is finite as well.

Together with the finite initial entropy, these verify the hypotheses
of the global density and logarithmic-gradient estimate in
\cite[Theorem~7.4.1]{Kri1}. It gives
$P_t(dx)=u_t^{\varepsilon,\delta}(x)\,dx$ for almost every $t$ and
\[
\int_0^T\int_{\R^d}
\frac{|\nabla u^{\varepsilon,\delta}|^2}
{u^{\varepsilon,\delta}}\,dx\,dt<\infty,
\]
where the quotient is set to zero on the zero set of the density.
The entropy estimate in the proof of the same theorem gives an
essential upper bound on
$\int u_t^{\varepsilon,\delta}\log u_t^{\varepsilon,\delta}\,dx$.
Together with the second-moment bound, this gives a uniform bound on
the relative entropy with respect to the standard Gaussian measure.
Lower semicontinuity of relative entropy and narrow continuity of
$P_t$ show that $P_t$ has a density also at the remaining times:
approximate each such time by times for which the entropy bound holds.
This is the narrowly continuous representative used above.

The identity $|\nabla\sqrt{u}|^2=|\nabla u|^2/(4u)$ and unit mass
give the claimed Sobolev regularity. Substituting the densities into
the weak equation gives Definition~\ref{weak nonlocal}, with the
same initial condition. The regularity constants may depend on
$\varepsilon$ and $\delta$.
\end{proof}

\section{Nonlocal-to-local limit}\label{4}
\subsection{Convergence as $\varepsilon\to0$}

The goal of this section is to prove that, as $\varepsilon,\delta\to0$,
weak solutions of \eqref{nonlocal},
\[
\begin{aligned}
\partial_t u^{\varepsilon,\delta}
&=\operatorname{div}\bigl(u^{\varepsilon,\delta}
(\nabla V+\nabla W*u^{\varepsilon,\delta})\bigr)\\
&\quad+\Delta\bigl(u^{\varepsilon,\delta}
\varphi_\varepsilon*(\varphi_\varepsilon*u^{\varepsilon,\delta})^{m-1}
+\delta u^{\varepsilon,\delta}\bigr),
\end{aligned}
\]
converge, along subsequences, to weak solutions of the granular medium
equation \eqref{granular medium equation},
\[
\partial_t u=\operatorname{div}\bigl(u(\nabla V+\nabla W*u)\bigr)
+\Delta u^m.
\]
As mentioned in the introduction, we first let $\varepsilon\to0$ with $\delta$ fixed,
and then let $\delta\to0$.

The following consequence of the assumptions will be useful:

\begin{lemma}
   Assume $(V)$ and $(W)$. Then
\begin{equation}\label{eq:limit-gradient-radial}
\begin{aligned}
|\nabla V(x)|^{\frac{2q}{2q-1}}
&\le C\bigl[x\cdot\nabla V(x)+c_V(1+|x|^2)\bigr],\\
|\nabla W(x)|^{\frac{2q}{2q-1}}
&\le C\bigl[x\cdot\nabla W(x)+c_W(1+|x|^2)\bigr].
\end{aligned}
\end{equation}

\end{lemma}
\begin{proof}   
Let $F$ stand for $V$ or $W$. If $\nabla F(x)\ne0$, apply inequality
\[
(\nabla F(x)-\nabla F(y))\cdot(x-y)
\geq -K_F|x-y|^2-C_{2,F}.
\]
to $x$ and
$y=r\nabla F(x)/|\nabla F(x)|$, where $r\geq1$. This gives
\[
r|\nabla F(x)|\leq x\cdot\nabla F(x)
  +|\nabla F(y)|\,|x-y|+K_F|x-y|^2+C_{2,F}.
\]
By growth assumptions and Young's inequality we have
\[
|\nabla F(y)||x-y|\leq Cr^q(|x|+r)\leq\frac{1}{2}|x|^2+Cr^{2q}.
\]
Similarly
\[K_F|x-y|^2\leq \left(K_F+\frac{1}{2}\right)|x|^2+Cr^{2q}.\]
These estimates together give us
\begin{align*}
r|\nabla F(x)|&\leq x\cdot\nabla F(x)+(K_F+1)|x|^2+Cr^{2q}\\
&\leq x\cdot\nabla F(x)+c_F(1+|x|^2)+Cr^{2q}.
\end{align*}
Substituting $r=
\bigl[x\cdot\nabla F(x)+c_F(1+|x|^2)\bigr]^{1/(2q)}$ we obtain the thesis.
\end{proof}

We begin with the regularity of the nonlocal solutions. The following
lemma follows from the estimates in \cite[Theorem~7.4.1]{Kri1},
as verified at the end of Section~\ref{3}.

\begin{lemma}\label{limit:lem:fixed-regularity}
For fixed $\varepsilon,\delta>0$, the chosen nonlocal solution satisfies
\[
\sqrt{u^{\varepsilon,\delta}}\in L^2(0,T;H^1(\R^d)),
\qquad
u^{\varepsilon,\delta}\in L^2(0,T;W^{1,1}(\R^d)).
\]
In particular, $u_t^{\varepsilon,\delta}\in W^{1,1}(\R^d)$ for
almost every $t$.
\end{lemma}
\begin{proof}
The first assertion is Corollary~\ref{chaos:cor:nonlocal-existence}.
The Sobolev chain rule and Cauchy--Schwarz give
\begin{equation}\label{eq:limit-gradient-density}
\begin{aligned}
\|\nabla u_t^{\varepsilon,\delta}\|_{L^1}
&=2\int_{\R^d}\sqrt{u_t^{\varepsilon,\delta}}\,
|\nabla\sqrt{u_t^{\varepsilon,\delta}}|\,dx\\
&\le2\|\nabla\sqrt{u_t^{\varepsilon,\delta}}\|_{L^2}.
\end{aligned}
\end{equation}
Square this inequality and integrate in time. The $L^1$ norm of the
density is one. The constants here need not be uniform in
$\varepsilon,\delta$.
\end{proof}

Before the main estimates, we recall two elementary facts which will
be useful in this section.

\begin{lemma}\label{limit:lem:entropy-lower}
For every probability density $\rho$ with finite second moment,
\begin{equation}\label{eq:limit-entropy-lower}
H(\rho)\ge-\frac d2\log(2\pi)-\frac12\int_{\R^d}|x|^2\rho(x)\,dx,
\end{equation}
where the entropy is allowed to be $+\infty$.
\end{lemma}
\begin{proof}
Consider the standard Gaussian density
\[
\mathcal M(x)=(2\pi)^{-d/2}\exp(-|x|^2/2).
\]
Jensen's inequality for $r\log r$ gives
\[
\begin{aligned}
\int_{\R^d}\rho\log\frac{\rho}{\mathcal M}\,dx
&=\int_{\R^d}\frac{\rho}{\mathcal M}
\log\frac{\rho}{\mathcal M}\,\mathcal M\,dx\\
&\ge
\left(\int_{\R^d}\frac{\rho}{\mathcal M}\mathcal M\,dx\right)
\log\left(\int_{\R^d}\frac{\rho}{\mathcal M}\mathcal M\,dx\right)=0.
\end{aligned}
\]
The negative part of $\rho\log\rho$ is integrable under the
second-moment assumption. Substituting $\log\mathcal M$ yields
\eqref{eq:limit-entropy-lower}.
\end{proof}

\begin{lemma}\label{limit:lem:interpolation}
Let $\Omega\subseteq\R^d$ be open and $1<m<p<\infty$. If
$f\in L^1(\Omega)\cap L^p(\Omega)$, then
\[
\|f\|_{L^m(\Omega)}
\le\|f\|_{L^1(\Omega)}^\theta
\|f\|_{L^p(\Omega)}^{1-\theta},
\qquad
\frac1m=\theta+\frac{1-\theta}{p}.
\]
Moreover, if $\rho$ is a smooth probability density on $\R^d$ and
$\rho^{m/2}\in H^1(\R^d)$, then
\begin{align}
\int_{\R^d}\rho^{m+2/d}\,dx
&\le C\|\nabla\rho^{m/2}\|_{L^2}^2,
\label{eq:limit-higher-integrability}\\
\int_{\R^d}\rho^m\,dx
&\le C\|\nabla\rho^{m/2}\|_{L^2}^{
\frac{2d(m-1)}{d(m-1)+2}},
\label{eq:limit-nonlinear-interpolation}
\end{align}
where constants depend on dimension and $m$.
These estimates hold in every dimension $d\ge1$.
\end{lemma}
\begin{proof}
Apply H\"older's inequality with conjugate exponents
$1/(m\theta)$ and $p/(m(1-\theta))$:
\[
\int_\Omega|f|^m
\le
\left(\int_\Omega|f|\right)^{m\theta}
\left(\int_\Omega|f|^p\right)^{m(1-\theta)/p}.
\]
Taking the $m$th root proves the interpolation inequality.

For $d\ge3$, H\"older's inequality, unit mass, and
Theorem~\ref{thm:prelim-gns} give
\[
\begin{aligned}
\int_{\R^d}\rho^{m+\frac{2}{d}}\,dx
&\le
\left(\int_{\R^d}\rho^{md/(d-2)}\,dx\right)^{(d-2)/d}
\left(\int_{\R^d}\rho\,dx\right)^{2/d}\\
&\le C\|\nabla\rho^{m/2}\|_{L^2}^2.
\end{aligned}
\]
For $d=2$, use the embedding $W^{1,1}(\R^2)\subset L^2(\R^2)$ and Cauchy-Schwarz inequality:
\[
\|\rho^{(m+1)/2}\|_{L^2}
\le C\|\nabla\rho^{(m+1)/2}\|_{L^1}
\le C\|\rho^{1/2}\|_{L^2}\|\nabla\rho^{m/2}\|_{L^2}.
\]
For $d=1$, the same argument with the embedding into $L^\infty$
gives
\[
\int_{\R}\rho^{m+2}\,dx
\le\|\rho^{(m+1)/2}\|_\infty^2\int_{\R}\rho\,dx
\le C\|\nabla\rho^{m/2}\|_{L^2}^2.
\]
This proves \eqref{eq:limit-higher-integrability}.
Interpolating between $L^1$ and $L^{m+2/d}$ now gives
\[
\int_{\R^d}\rho^m\,dx
\le\left(\int_{\R^d}\rho^{m+2/d}\,dx\right)^{
\frac{d(m-1)}{d(m-1)+2}},
\]
which proves \eqref{eq:limit-nonlinear-interpolation}.
\end{proof}

Now we can prove the crucial estimate of this section.

\begin{lemma}\label{limit:lem:entropy-moments}
Assume $(V),(W),(\mu_0)$ and \eqref{joint}. Then there
is a constant $C$ independent of $0<\varepsilon,\delta\leq1$ such that
\begin{align}
&\sup_{t\in[0,T]}\left[
H(u_t^{\varepsilon,\delta})
+\int_{\R^d}|x|^2u_t^{\varepsilon,\delta}\,dx\right]\le C,
\label{eq:limit-entropy-moment-bound}\\
&\sup_{t\in[0,T]}\int_{\R^d}|x|^2u_t^{\varepsilon,\delta}\,dx
+\|(\varphi_\varepsilon*u^{\varepsilon,\delta})^{m/2}
\|_{L^2(0,T;H^1)}^2
+\delta\|\sqrt{u^{\varepsilon,\delta}}\|_{L^2(0,T;H^1)}^2
\le C.
\label{eq:limit-sobolev-bound}
\end{align}
In addition,
\begin{equation}\label{eq:limit-radial-bound}
\begin{aligned}
&\int_0^T\int_{\R^d}
\bigl[x\cdot\nabla V(x)+c_V(1+|x|^2)\bigr]
u_t^{\varepsilon,\delta}(x)\,dx\,dt\\
&\quad+\frac12\int_0^T\iint_{\R^d\times\R^d}
\bigl[(x-y)\cdot\nabla W(x-y)+c_W(1+|x-y|^2)\bigr]\\
&\hspace{47mm}\cdot
u_t^{\varepsilon,\delta}(x)u_t^{\varepsilon,\delta}(y)
\,dx\,dy\,dt\le C.
\end{aligned}
\end{equation}
\end{lemma}

\begin{proof}
Since the initial entropy and second moment are finite, we use the
tests $\log u^{\varepsilon,\delta}+1$ and $|x|^2$.
Their admissibility, by truncating the logarithm and the spatial
variable, is proved in Appendix~\ref{app:test-functions}.
We split the calculation into two parts.

The entropy calculation gives, in integrated form,
\begin{equation}
\begin{aligned}
H(\u_t)-H(u_0)=-\int_0^t\int_{\R^d}(\nabla V+\nabla W*\u)\u\frac{\nabla\u}{\u}\\-\int_0^t\delta\int_{\R^d}\frac{|\nabla\u|^2}{\u}-\int_0^t\int_{\R^d}\frac{|\nabla \u|^2}{\u}(\varphi_\epsilon*(\varphi_\epsilon*\u)^{m-1})\\-\int_0^t\int_{\R^d}\frac{\nabla\u}{\u}\u\nabla\Big(\varphi_\epsilon*(\varphi_\epsilon*\u)^{m-1}\Big)
\end{aligned}
\end{equation}
Using moment bounds of $\u_t$ for fixed $\epsilon,\delta$ we can integrate by parts in term involving $V$ and $W$. Moreover, moving the convolution onto
$\nabla u^{\varepsilon,\delta}$:
\[
\begin{aligned}
&\int_{\R^d}\nabla u^{\varepsilon,\delta}\cdot
\nabla\Big(\varphi_\varepsilon*
(\varphi_\varepsilon*u^{\varepsilon,\delta})^{m-1}\Big)\,\text{d}x\\
&\quad=(m-1)\int_{\R^d}\nabla\u \Big(\varphi_\epsilon*(\varphi_\epsilon*\u)^{m-2}\nabla(\varphi_\epsilon*\u)\Big)\,\text{d}x\\
&\quad=(m-1)\int_{\R^d}
(\varphi_\varepsilon*u^{\varepsilon,\delta})^{m-2}
|\nabla(\varphi_\varepsilon*u^{\varepsilon,\delta})|^2\,\text{d}x\\
&\quad=\frac{4(m-1)}{m^2}
\|\nabla(\varphi_\varepsilon*u^{\varepsilon,\delta})^{m/2}\|_{L^2}^2.
\end{aligned}
\]
To estimate the Laplacians, use unit mass and the joint condition:
\[
\begin{aligned}
&\int_{\R^d}(\Delta V+\Delta W*u^{\varepsilon,\delta})
u^{\varepsilon,\delta}\,dx\\
&=\frac12\iint_{\R^d\times\R^d}
\bigl[\Delta V(x)+\Delta V(y)+2\Delta W(x-y)\bigr]
u^{\varepsilon,\delta}(x)u^{\varepsilon,\delta}(y)\,dx\,dy\\
&\le A\int_{\R^d}
\bigl[x\cdot\nabla V(x)+c_V(1+|x|^2)\bigr]
u^{\varepsilon,\delta}(x)\,dx\\
&\quad+\frac A2\iint_{\R^d\times\R^d}
\bigl[(x-y)\cdot\nabla W(x-y)+c_W(1+|x-y|^2)\bigr]\\
&\hspace{49mm}\cdot u^{\varepsilon,\delta}(x)
u^{\varepsilon,\delta}(y)\,dx\,dy.
\end{aligned}
\]
All terms are integrable for fixed $\varepsilon,\delta$ by the
moment bounds from Section~\ref{3}.

Testing by $|x|^2$ and symmetrizing the interaction term gives
\begin{equation}\label{eq:limit-second-moment}
\begin{aligned}
\frac{d}{dt}\int_{\R^d}|x|^2u^{\varepsilon,\delta}\,dx
&=-2\int_{\R^d}x\cdot\nabla V(x)u^{\varepsilon,\delta}(x)\,dx\\
&\quad-\iint_{\R^d\times\R^d}(x-y)\cdot\nabla W(x-y)
u^{\varepsilon,\delta}(x)u^{\varepsilon,\delta}(y)\,dx\,dy\\
&\quad+2d\delta
+2d\int_{\R^d}(\varphi_\varepsilon*u^{\varepsilon,\delta})^m\,dx
\end{aligned}
\end{equation}
for almost every $t$. Since
\[
\iint|x-y|^2u^{\varepsilon,\delta}(x)
u^{\varepsilon,\delta}(y)\,dx\,dy
\le2\int|x|^2u^{\varepsilon,\delta}(x)\,dx,
\]
adding to \eqref{eq:limit-second-moment} the quadratic corrections in such way that the following formula match the Laplacian estimates we obtained
 yields
\[
\begin{aligned}
&\frac{d}{dt}\int_{\R^d}|x|^2u^{\varepsilon,\delta}\,dx
+2\int_{\R^d}\bigl[x\cdot\nabla V(x)+c_V(1+|x|^2)\bigr]
u^{\varepsilon,\delta}(x)\,dx\\
&\quad+\iint_{\R^d\times\R^d}
\bigl[(x-y)\cdot\nabla W(x-y)+c_W(1+|x-y|^2)\bigr]
u^{\varepsilon,\delta}(x)u^{\varepsilon,\delta}(y)\,dx\,dy\\
&\le C\left(1+\int_{\R^d}|x|^2u^{\varepsilon,\delta}\,dx\right)
+2d\int_{\R^d}(\varphi_\varepsilon*u^{\varepsilon,\delta})^m\,dx.
\end{aligned}
\]

Choose $\gamma=A+1$. Multiply the last inequality by $\gamma$, integrate with respect to time
and add it to the entropy inequality.
By \eqref{eq:limit-nonlinear-interpolation} and Young's inequality, since the exponent in the first line is strictly smaller than two, we have
\[
\begin{aligned}
2d\gamma\int_{\R^d}(\varphi_\varepsilon*u^{\varepsilon,\delta})^m\,dx
&\le C\|\nabla(\varphi_\varepsilon*u^{\varepsilon,\delta})^{m/2}
\|_{L^2}^{\frac{2d(m-1)}{d(m-1)+2}}\\
&\le C+\frac{2(m-1)}{m^2}
\|\nabla(\varphi_\varepsilon*u^{\varepsilon,\delta})^{m/2}\|_{L^2}^2.
\end{aligned}
\]
Thus,
\begin{equation}\label{eq:limit-combined-estimate}
\begin{aligned}
&H(u_t^{\varepsilon,\delta})
+\gamma\int_{\R^d}|x|^2u_t^{\varepsilon,\delta}\,dx\\
&\quad+\frac{2(m-1)}{m^2}\int_0^t
\|\nabla(\varphi_\varepsilon*u^{\varepsilon,\delta})^{m/2}\|_{L^2}^2\,dr
+4\delta\int_0^t\|\nabla\sqrt{u^{\varepsilon,\delta}}\|_{L^2}^2\,dr\\
&\quad+(2\gamma-A)\int_0^t\int_{\R^d}
\bigl[x\cdot\nabla V(x)+c_V(1+|x|^2)\bigr]
u^{\varepsilon,\delta}(x)\,dx\,dr\\
&\quad+\frac{2\gamma-A}{2}\int_0^t\iint_{\R^d\times\R^d}
\bigl[(x-y)\cdot\nabla W(x-y)+c_W(1+|x-y|^2)\bigr]\\
&\hspace{51mm}\cdot
u^{\varepsilon,\delta}(x)u^{\varepsilon,\delta}(y)\,dx\,dy\,dr\\
&\le H(u_0)+\gamma\int_{\R^d}|x|^2u_0\,dx
+C\int_0^t\left(1+\int_{\R^d}|x|^2u_r^{\varepsilon,\delta}\,dx\right)\,dr.
\end{aligned}
\end{equation}
By \eqref{eq:limit-entropy-lower}, the first line is bounded from below up to an additive constant, since $\gamma\geq\frac{1}{2}$. Moreover, $2\gamma-A>0$ and we can use Gronwall's inequality to obtain uniform second-moment bound.
Returning to \eqref{eq:limit-combined-estimate} gives the entropy,
gradient, and \eqref{eq:limit-radial-bound} bounds. Finally,
\eqref{eq:limit-nonlinear-interpolation} gives the time-integrated
$L^m$ bound for $\varphi_\varepsilon*u^{\varepsilon,\delta}$.
Unit mass supplies the $L^2$ part of the norm of
$\sqrt{u^{\varepsilon,\delta}}$.
The entropy bounds extend to every time by lower semicontinuity
of relative entropy and narrow continuity.
\end{proof}

An important consequence, used in both limiting procedures, is
\begin{equation}\label{eq:limit-drift-power}
\begin{aligned}
&\int_0^T\int_{\R^d}
|\nabla V(x)|^{\frac{2q}{2q-1}}u^{\varepsilon,\delta}(x)\,dx\,dt\\
&\quad+\int_0^T\iint_{\R^d\times\R^d}
|\nabla W(x-y)|^{\frac{2q}{2q-1}}
u^{\varepsilon,\delta}(x)u^{\varepsilon,\delta}(y)\,dx\,dy\,dt
\le C.
\end{aligned}
\end{equation}
Indeed, apply \eqref{eq:limit-gradient-radial} to
\eqref{eq:limit-radial-bound}. Jensen's inequality then gives the
same bound for
$|\nabla V+\nabla W*u^{\varepsilon,\delta}|^{2q/(2q-1)}$
integrated against $u^{\varepsilon,\delta}\,dx\,dt$.

Using these estimates, we can prove strong compactness.
We recall the version of the Aubin--Lions compactness lemma from \cite[Corollary~4]{Sim1}.

\begin{lemma}[Aubin--Lions compactness lemma]\label{limit:lem:Aubin-Lions}
Let $X_0,X,X_1$ be Banach spaces with
$X_0\subset\subset X\hookrightarrow X_1$, and let $1\le p<\infty$.
A family bounded in $L^p(0,T;X_0)$ is relatively compact in
$L^p(0,T;X)$ if its distributional time derivatives are bounded
in $L^1(0,T;X_1)$.
\end{lemma}

\begin{theorem}\label{limit:thm:epsilon}
Assume $(V),(W),(\mu_0)$ and \eqref{joint}, fix
$\delta\in(0,1]$. Every sequence $\varepsilon\to0$ has a subsequence
such that
\[
u^{\varepsilon,\delta}\to u^\delta
\quad\text{in }L^1((0,T)\times\R^d),
\qquad
\varphi_\varepsilon*u^{\varepsilon,\delta}\to u^\delta
\quad\text{in }L^m((0,T)\times\R^d).
\]
The corresponding measure curves converge narrowly, uniformly on
$[0,T]$, and $u^\delta$ is a weak solution of \eqref{delta} with
initial datum $u_0$.
\end{theorem}

\begin{proof}
The proof has two main parts. First, we show strong compactness in
order to pass to the limit in the nonlocal diffusion term. Then we
show that the limit satisfies \eqref{delta}.
We work first on a ball $B_R$, using
\[
W^{1,1}(B_R)\subset\subset L^1(B_R)
\hookrightarrow (W^{2,\infty}_0(B_R))^*.
\]
Here $W^{2,\infty}_0(B_R)$ means the closure of $C_c^\infty(B_R)$
in $W^{2,\infty}$; its elements have admissible zero extensions.

\begin{enumerate}
\item \textit{Strong compactness.}
By \eqref{eq:limit-gradient-density} and
Lemma~\ref{limit:lem:entropy-moments},
\[
\int_0^T\|\nabla u^{\varepsilon,\delta}\|_{L^1}^2\,dt
\le4\int_0^T\|\nabla\sqrt{u^{\varepsilon,\delta}}\|_{L^2}^2\,dt
\le\frac C\delta.
\]
Thus $u^{\varepsilon,\delta}$ and
$\varphi_\varepsilon*u^{\varepsilon,\delta}$ are bounded in
$L^2(0,T;W^{1,1}(B_R))$ for fixed $\delta$.

Take $\psi\in W^{2,\infty}_0(B_R)$ with norm at most one and test the equation by $\psi*\varphi_\epsilon$. By evenness of the kernel,
\[
\begin{aligned}
&\left\langle\partial_t(\varphi_\varepsilon*
u^{\varepsilon,\delta}),\psi\right\rangle\\
&=\int_{\R^d}
\left[\delta u^{\varepsilon,\delta}
+u^{\varepsilon,\delta}\varphi_\varepsilon*
(\varphi_\varepsilon*u^{\varepsilon,\delta})^{m-1}\right]
(\varphi_\varepsilon*\Delta\psi)\,dx\\
&\quad-\int_{\R^d}u^{\varepsilon,\delta}
(\nabla V+\nabla W*u^{\varepsilon,\delta})\cdot
(\varphi_\varepsilon*\nabla\psi)\,dx.
\end{aligned}
\]
The absolute value is bounded by
\[
\delta+\int_{\R^d}(\varphi_\varepsilon*u^{\varepsilon,\delta})^m\,dx
+\int_{\R^d}u^{\varepsilon,\delta}
|\nabla V+\nabla W*u^{\varepsilon,\delta}|\,dx.
\]
Its time integral is bounded independently of $\varepsilon$ by
\eqref{eq:limit-sobolev-bound} and \eqref{eq:limit-drift-power}.
This proves the required bound for the time derivative.
The same calculation without convolution gives a bound for
$\partial_tu^{\varepsilon,\delta}$.

Lemma~\ref{limit:lem:Aubin-Lions} gives, along a subsequence,
strong convergence of $\varphi_\varepsilon*u^{\varepsilon,\delta}$
in $L^1((0,T)\times B_R)$.
The translation estimate for $W^{1,1}(\R^d)$ gives
\[
\begin{aligned}
&\|u^{\varepsilon,\delta}
-\varphi_\varepsilon*u^{\varepsilon,\delta}\|_
{L^1((0,T)\times\R^d)}\\
&\quad\le
\varepsilon C
\int_0^T\|\nabla u^{\varepsilon,\delta}\|_{L^1}\,dt
\le C_\delta\varepsilon.
\end{aligned}
\]
Consequently, the density itself has the same local $L^1$ limit,
which we denote by $u^\delta$.

Next, \eqref{eq:limit-higher-integrability} gives
\[
\int_0^T\int_{\R^d}
(\varphi_\varepsilon*u^{\varepsilon,\delta})^{m+2/d}\,dx\,dt
\le C.
\]
The interpolation inequality of
Lemma~\ref{limit:lem:interpolation} upgrades local $L^1$ convergence
to local $L^m$ convergence of the convolved densities:
\begin{align*}
||\u*\varphi_\epsilon-u^\delta||_{L^m((0,T)\times B_R)}\leq||\u*\varphi_\epsilon-u^\delta||^\theta_{L^1((0,T)\times B_R)}\\ \times ||\u*\varphi_\epsilon-u^\delta||^{1-\theta}_{L^{m+2/d}((0,T)\times B_R)}.
\end{align*}

\item \textit{Convergence to a solution of \eqref{delta}.}
Take increasing balls $B_1,B_2,\ldots$ with radii tending to infinity.
Successively extract subsequences converging on these balls and
then take the diagonal subsequence. We keep the same notation for it.
The second-moment bound and local convergence give us global convergence.

The time compactness argument in
Proposition~\ref{limit:prop:time-continuity} below uses only
Lemma~\ref{limit:lem:entropy-moments}, \eqref{eq:limit-drift-power},
and the weak equation. It therefore applies at this stage and,
after a further extraction, gives uniform narrow convergence of
the measure curves. This identifies $u^\delta$ at every time,
preserves unit mass, and gives the initial and terminal terms
in the weak formulation. The uniform entropy and moment bounds
ensure that these time marginals are densities.

For $\psi\in C_c^\infty([0,T]\times\R^d)$, the nonlocal equation reads
\[
\begin{aligned}
&-\int_0^T\int_{\R^d}u^{\varepsilon,\delta}\partial_t\psi\,dx\,dt
+\int_{\R^d}u_T^{\varepsilon,\delta}\psi_T\,dx
-\int_{\R^d}u_0\psi_0\,dx\\
&=-\int_0^T\int_{\R^d}u^{\varepsilon,\delta}
(\nabla V+\nabla W*u^{\varepsilon,\delta})\cdot\nabla\psi\,dx\,dt\\
&\quad+\int_0^T\int_{\R^d}u^{\varepsilon,\delta}
\bigl[\varphi_\varepsilon*
(\varphi_\varepsilon*u^{\varepsilon,\delta})^{m-1}\bigr]
\Delta\psi\,dx\,dt\\
&\quad+\delta\int_0^T\int_{\R^d}u^{\varepsilon,\delta}\Delta\psi\,dx\,dt.
\end{aligned}
\]
The time, linear diffusion, and external-potential terms converge
by the convergence already obtained.

For the nonlinear diffusion term, use the symmetry of the kernel
and add and subtract the same expression:
\[
\begin{aligned}
&\int_0^T\int_{\R^d}u^{\varepsilon,\delta}
\bigl[\varphi_\varepsilon*
(\varphi_\varepsilon*u^{\varepsilon,\delta})^{m-1}\bigr]
\Delta\psi\,dx\,dt\\
&=\int_0^T\int_{\R^d}
(\varphi_\varepsilon*u^{\varepsilon,\delta})^m\Delta\psi\,dx\,dt\\
&\quad+\int_0^T\int_{\R^d}
\bigl[\varphi_\varepsilon*(u^{\varepsilon,\delta}\Delta\psi)
-(\varphi_\varepsilon*u^{\varepsilon,\delta})\Delta\psi\bigr]
(\varphi_\varepsilon*u^{\varepsilon,\delta})^{m-1}\,dx\,dt.
\end{aligned}
\]
The first term converges to $\int_0^T\int(u^\delta)^m\Delta\psi$.
Since $\varphi$ has compact support and $\Delta\psi$ is Lipschitz,
\begin{align*}
\left|\varphi_\varepsilon*(u^{\varepsilon,\delta}\Delta\psi)
-(\varphi_\varepsilon*u^{\varepsilon,\delta})\Delta\psi\right|\\=\leq
\int_{\mathbb{R}^d} \varphi_\epsilon(x-y) u_\epsilon(t, y) \big| \Delta \psi(t, y) - \Delta \psi(t, x) \big| \, \text{d}y
\leq 
(\varphi_\varepsilon*u^{\varepsilon,\delta})\,C_\psi\varepsilon.
\end{align*}
The absolute value of the second integral is therefore at most
\[
C_\psi\varepsilon\int_0^T\int_{\R^d}
(\varphi_\varepsilon*u^{\varepsilon,\delta})^m\,dx\,dt,
\]
which tends to zero.

It remains to treat the interaction term. Let $\chi_R$ be a smooth
cutoff equal to one on $B_R$ and zero outside $B_{2R}$, and replace
$\nabla W(x-y)$ by $\chi_R(x-y)\nabla W(x-y)$.
For fixed $R$, this kernel is bounded and continuous.
The interaction term is an integral against
$u_t^{\varepsilon,\delta}(x)u_t^{\varepsilon,\delta}(y)\,dx\,dy\,dt$,
so its truncated version converges by the strong $L^1$ convergence.
For the discarded part, H\"older's inequality gives
\[
\begin{aligned}
&\int_0^T\iint_{\{|x-y|>R\}}
|\nabla W(x-y)|u_t^{\varepsilon,\delta}(x)
u_t^{\varepsilon,\delta}(y)\,dx\,dy\,dt\\
&\le
\left(\int_0^T\iint
|\nabla W(x-y)|^{\frac{2q}{2q-1}}
u_t^{\varepsilon,\delta}(x)u_t^{\varepsilon,\delta}(y)
\,dx\,dy\,dt\right)^{\frac{2q-1}{2q}}\\
&\quad\cdot
\left(\int_0^T\iint_{\{|x-y|>R\}}
u_t^{\varepsilon,\delta}(x)u_t^{\varepsilon,\delta}(y)
\,dx\,dy\,dt\right)^{1/(2q)}
\le C R^{-1/q}.
\end{aligned}
\]
The same bound holds for the limit, by lower semicontinuity.
Letting first $\varepsilon\to0$ and then $R\to\infty$ proves
convergence of the full interaction term.

We have obtained \eqref{delta}. The drift is globally integrable
by the limiting version of \eqref{eq:limit-drift-power}, and
$(u^\delta)^m$ is integrable by strong $L^m$ convergence.
Thus all requirements of Definition~\ref{weak delta} are satisfied.
\end{enumerate}
\end{proof}

Without a uniqueness result for \eqref{delta}, the conclusion is
subsequential convergence. If the weak solution is unique in the
class above, compactness identifies every subsequential limit and
the whole family converges.

\subsection{Convergence as $\delta\to0$}

Having obtained solutions of \eqref{delta}, we now let $\delta$
tend to zero. We begin with the estimates needed for compactness.

\begin{lemma}\label{limit:lem:delta-bounds}
The solutions $u^\delta$ constructed in
Theorem~\ref{limit:thm:epsilon} satisfy
\[
\sup_{t\in[0,T]}\int_{\R^d}|x|^2u_t^\delta\,dx
+\|(u^\delta)^{m/2}\|_{L^2(0,T;H^1)}^2
+\delta\|\sqrt{u^\delta}\|_{L^2(0,T;H^1)}^2\le C,
\]
with $C$ independent of $1>\delta>0$.
The entropy and drift estimates
\eqref{eq:limit-entropy-moment-bound},
\eqref{eq:limit-radial-bound}, and \eqref{eq:limit-drift-power}
also hold with $u^\delta$ in place of $u^{\varepsilon,\delta}$.
Moreover,
\[
\int_0^T\int_{\R^d}(u^\delta)^{m+2/d}\,dx\,dt\le C.
\]
\end{lemma}
\begin{proof}
The strong $L^m$ convergence of the convolved densities implies
strong $L^2$ convergence of their powers $m/2$.
Weak lower semicontinuity in $L^2(0,T;H^1)$ gives the first gradient
bound. Similarly,
\[
\|\sqrt{u^{\varepsilon,\delta}}-\sqrt{u^\delta}\|_{L^2}^2
\le\|u^{\varepsilon,\delta}-u^\delta\|_{L^1}
\]
gives the bound for $\sqrt{u^\delta}$.
The corrected radial integrands are continuous and nonnegative;
their bounds pass to the limit by Fatou's lemma and convergence
of the product measures. The same reasoning applies to the
drift powers. The second moment and the relative entropy are
lower semicontinuous under narrow convergence.
Finally, apply \eqref{eq:limit-higher-integrability} to $\u$ and use lower semicontinuity to obtain uniform in $\delta$ integrability of $(m+2/d)$-power of $u^\delta$.
\end{proof}

We can now prove strong compactness of $u^\delta$ in $L^m$.
The proof is slightly different, since the bound for
$\nabla u^\delta$ obtained from $\nabla\sqrt{u^\delta}$ depends
on $\delta$. To do so, we use fractional Sobolev spaces.

\begin{theorem}\label{limit:thm:delta}
Assume $(V),(W),(\mu_0)$ and \eqref{joint}. Then
every sequence $\delta\to0$ has a subsequence such that
\[
u^\delta\to u
\quad\text{strongly in }L^m((0,T)\times\R^d)
\text{ and in }L^1((0,T)\times\R^d).
\]
The measure curves converge narrowly, uniformly on $[0,T]$.
The limit is a weak solution of \eqref{granular medium equation}
with initial datum $u_0$.
\end{theorem}

\begin{proof}
The proof follows that of Theorem~\ref{limit:thm:epsilon}.
The additional point is spatial compactness. We control
$(u^\delta)^{m/2}$ in $L^2(0,T;H^1)$, but we do not have a
corresponding estimate for its time derivative.
Instead, we use this spatial bound to estimate $u^\delta$ in
a fractional Sobolev space on every ball $B_R$.

We consider two cases.
\begin{enumerate}
\item \textit{Case $m\ge2$.}
The inequality $|A^\alpha-B^\alpha|\le|A-B|^\alpha$ for
$A,B\ge0$ and $0<\alpha\le1$ gives
\[
|u^\delta(x)-u^\delta(y)|^m
\le|(u^\delta(x))^{m/2}-(u^\delta(y))^{m/2}|^2.
\]
For $0<s<2/m$, it follows that
\[
\begin{aligned}
{}[u^\delta]_{W^{s,m}(B_R)}^m
&\le\iint_{B_R\times B_R}
\frac{|(u^\delta(x))^{m/2}-(u^\delta(y))^{m/2}|^2}
{|x-y|^{d+sm}}\,dx\,dy\\
&=[(u^\delta)^{m/2}]_{H^{sm/2}(B_R)}^2
\le C\|(u^\delta)^{m/2}\|_{H^1(B_R)}^2.
\end{aligned}
\]
The strict inequality $sm/2<1$ is essential for this use of the
fractional seminorm. Integrating in time gives a bound in
$L^m(0,T;W^{s,m}(B_R))$.

\item \textit{Case $1<m<2$.}
For $A,B\ge0$ and $\alpha>1$,
\[
|A^\alpha-B^\alpha|
\le C_\alpha(A^{\alpha-1}+B^{\alpha-1})|A-B|.
\]
Taking $\alpha=2/m$ gives
\[
\begin{aligned}
|u^\delta(x)-u^\delta(y)|
&\le C\bigl[(u^\delta(x))^{(2-m)/2}
+(u^\delta(y))^{(2-m)/2}\bigr]\\
&\hspace{14mm}\cdot
|(u^\delta(x))^{m/2}-(u^\delta(y))^{m/2}|.
\end{aligned}
\]
Fix $0<s<\sigma<1$. Cauchy--Schwarz yields
\[
\begin{aligned}
{}[u^\delta]_{W^{s,1}(B_R)}
&\le C\left(\iint_{B_R\times B_R}
\frac{(u^\delta(x))^{2-m}+(u^\delta(y))^{2-m}}
{|x-y|^{d-2\sigma+2s}}\,dx\,dy\right)^{1/2}\\
&\quad\cdot
\left(\iint_{B_R\times B_R}
\frac{|(u^\delta(x))^{m/2}-(u^\delta(y))^{m/2}|^2}
{|x-y|^{d+2\sigma}}\,dx\,dy\right)^{1/2}.
\end{aligned}
\]
The second factor is the $H^\sigma$ seminorm and is bounded by
the $H^1$ norm. Since $\sigma>s$, the first squared factor is at most
\[
C\int_{B_R}(u^\delta)^{2-m}\,dx
\le C|B_R|^{m-1}
\left(\int_{B_R}u^\delta\,dx\right)^{2-m}
\le C|B_R|^{m-1}.
\]
Consequently,
\[
[u^\delta]_{W^{s,1}(B_R)}
\le C\|(u^\delta)^{m/2}\|_{H^1(B_R)}.
\]
Unit mass and integration in time give a bound in
$L^2(0,T;W^{s,1}(B_R))$.
\end{enumerate}

In both cases the spatial space is compactly embedded in $L^1(B_R)$;
see Theorem~\ref{thm:prelim-fractional-compactness}.
For $\psi\in W^{2,\infty}_0(B_R)$, the weak equation gives
\[
\begin{aligned}
|\langle\partial_tu^\delta,\psi\rangle|
&\le\|\psi\|_{W^{2,\infty}}
\left[\int_{\R^d}(u^\delta)^m\,dx+\delta
+\int_{\R^d}u^\delta|\nabla V+\nabla W*u^\delta|\,dx\right].
\end{aligned}
\]
The bracket has a uniform $L^1(0,T)$ bound.
Lemma~\ref{limit:lem:Aubin-Lions} therefore gives a strongly
convergent subsequence in $L^1((0,T)\times B_R)$.
The diagonal argument, the second-moment tail estimate, and
interpolation with the uniform $L^{m+2/d}$ bound give the asserted
global $L^1$ and $L^m$ convergence.

As in Theorem~\ref{limit:thm:epsilon},
Proposition~\ref{limit:prop:time-continuity} supplies uniform narrow
convergence of the measure curves. In the weak formulation,
$(u^\delta)^m\to u^m$ in $L^1$, whereas
\[
\left|\delta\int_0^T\int_{\R^d}u^\delta\Delta\psi\,dx\,dt\right|
\le\delta T\|\Delta\psi\|_\infty\longrightarrow0.
\]
The external-potential term converges locally, and the interaction
term converges by the same cutoff and $CR^{-1/q}$ tail estimate
as before. The limiting drift estimate gives global integrability.
The initial condition is preserved by uniform narrow convergence.
This proves taht limit satisfies Definition~\ref{weak}.
\end{proof}

As before, convergence of the whole family follows if uniqueness
is known in this class. The compactness argument alone gives
existence and subsequential convergence.

\subsection{Martingale solution}

In this subsection we connect the PDE and SDE solutions.
We show that the weak solution of
\eqref{granular medium equation} gives a martingale solution
of \eqref{SDE}. We use the superposition principle,
Theorem~\ref{thm:prelim-superposition}, proved in \cite{Tre1}.
We first complete the time compactness argument used above.

\begin{proposition}\label{limit:prop:nonlocal-continuity}
For fixed $\varepsilon,\delta>0$, the chosen curve
$(u_t^{\varepsilon,\delta}(x)\,dx)_{t\in[0,T]}$ is narrowly continuous.
\end{proposition}
\begin{proof}
Section~\ref{3} gives a continuous process
$Z^{\varepsilon,\delta}$ with these marginal laws.
For $f\in C_b(\R^d)$, dominated convergence gives
\[
\lim_{s\to t}\E f(Z_s^{\varepsilon,\delta})
=\E f(Z_t^{\varepsilon,\delta}).
\]
This is narrow continuity.
\end{proof}

For the limiting curves we use the Arzel\`a--Ascoli theorem in the
space of probability measures.

\begin{theorem}[Arzel\`a--Ascoli theorem]\label{limit:thm:Ascoli}
Let $K$ be a compact metric space with distance $d$, and let
$u_n:[0,T]\to K$ be continuous curves.
If
\[
d(u_n(s),u_n(t))\le\omega(|t-s|)
\quad\text{for all }n,s,t,
\qquad \lim_{h\downarrow0}\omega(h)=0,
\]
then a subsequence converges uniformly on $[0,T]$ to a continuous curve.
\end{theorem}

Before applying this theorem, we need the following time-integrability
estimate.

\begin{lemma}\label{limit:lem:time-integrability}
Let $p=1+\frac{2}{d(m-1)}$. Then
\[
\int_0^T
\left(\int_{\R^d}(\varphi_\varepsilon*
u^{\varepsilon,\delta})^m\,dx\right)^pdt\le C,
\qquad
\int_0^T\left(\int_{\R^d}(u^\delta)^m\,dx\right)^pdt\le C,
\]
uniformly in $0<\varepsilon,\delta\le1$.
\end{lemma}
\begin{proof}
By \eqref{eq:limit-nonlinear-interpolation},
\[
\left(\int_{\R^d}(\varphi_\varepsilon*
u^{\varepsilon,\delta})^m\,dx\right)^p
\le C\|\nabla(\varphi_\varepsilon*
u^{\varepsilon,\delta})^{m/2}\|_{L^2}^2,
\]
because $\frac{pd(m-1)}{d(m-1)+2}=1$.
Integrate in time and use Lemma~\ref{limit:lem:entropy-moments}.
The argument for $u^\delta$ is identical, using
Lemma~\ref{limit:lem:delta-bounds}.
\end{proof}

\begin{proposition}\label{limit:prop:time-continuity}
The nonlocal curves satisfy
\[
d_{BL}(u_t^{\varepsilon,\delta},u_s^{\varepsilon,\delta})
\le C|t-s|^\alpha,
\qquad
\alpha=\min\left\{\frac1{2q},\frac1{d(m-1)+2}\right\},
\]
with $C$ independent of $0<\varepsilon,\delta\le1$.
They are relatively compact for uniform narrow convergence on $[0,T]$.
The same statements hold for the curves $u^\delta$ constructed above.
Every limiting curve $u^\delta$ and $u$ is narrowly continuous.
\end{proposition}
\begin{proof}
The uniform second-moment bound implies that all time marginals
belong to a single compact subset of $\mathcal P(\R^d)$.
We prove equicontinuity using the bounded Lipschitz distance
\[
d_{BL}(\mu,\nu)=
\sup_{\substack{\|\psi\|_\infty\le1\\
\operatorname{Lip}(\psi)\le1}}
\left|\int_{\R^d}\psi\,d(\mu-\nu)\right|.
\]
Take an admissible $\psi$, a smooth nonnegative even mollifier
$\theta$ on $\R^d$, and $0<h\le1$.
Standard convolution estimates give
\[
\|\psi-\psi*\theta_h\|_\infty\le Ch,\qquad
\|\nabla(\psi*\theta_h)\|_\infty\le1,\qquad
\|\Delta(\psi*\theta_h)\|_\infty\le C/h.
\]
Let $\chi_R$ be a smooth cutoff equal to one on $B_R$ and zero
outside $B_{2R}$, with
$|\nabla\chi_R|\le C/R$ and $|\Delta\chi_R|\le C/R^2$.
For $\psi_{R,h}=(\psi*\theta_h)\chi_R$ and $R\ge1$,
\[
\|\nabla\psi_{R,h}\|_\infty\le C,\qquad
\|\Delta\psi_{R,h}\|_\infty\le C/h.
\]

We split the difference into three terms:
\[
\begin{aligned}
\left|\int_{\R^d}\psi(u_t^{\varepsilon,\delta}
-u_s^{\varepsilon,\delta})\,\text{d}x\right|
\leq \left|\int_{\R^d}(\psi-\psi*\theta_h)(\u_t-\u_s)\,\text{d}x\right|\\
+\left|\int_{\R^d}\psi*\theta_h(1-\chi_R)(\u_t-\u_s)\,\text{d}x\right|
+
\left|\int_{\R^d}\psi_{R,h}(u_t^{\varepsilon,\delta}
-u_s^{\varepsilon,\delta})\,dx\right|\\\leq
2Ch+\frac C{R^2}+
\left|\int\psi_{R,h}(u_t^{\varepsilon,\delta}
-u_s^{\varepsilon,\delta})\,dx\right|
\end{aligned}
\]
The first term is the mollification error; the second is the
spatial tail, estimated by the second moment.
For $s<t$, the weak equation gives
\[
\begin{aligned}
\left|\int\psi_{R,h}(u_t^{\varepsilon,\delta}
-u_s^{\varepsilon,\delta})\,\text{d}x\right|&\leq C\int_s^t\int u^{\varepsilon,\delta}
|\nabla V+\nabla W*u^{\varepsilon,\delta}|\,\text{d}x\,\text{d}r\\
&\quad+\frac C h\int_s^t
\left[\delta+\int(\varphi_\varepsilon*u^{\varepsilon,\delta})^m\,\text{d}x
\right]\,\text{d}r\\
&\le C|t-s|^{1/(2q)}
+\frac C h\bigl(|t-s|+|t-s|^{1-1/p}\bigr).
\end{aligned}
\]
Here the drift is estimated by H\"older's inequality and
\eqref{eq:limit-drift-power}:
\begin{align*}
    \int_s^t\int_{\R^d}u^{\varepsilon,\delta}
|\nabla V+\nabla W*u^{\varepsilon,\delta}|\,\text{d}x\,\text{d}r\\\leq\left(
\int_s^t\int_{\mathbb R^d}
\left|\nabla V+\nabla W*\u_r\right|^{\frac{2q}{2q-1}}
\u_r\,\text{d}x\,\text{d}r
\right)^{\frac{2q-1}{2q}}\cdot\left(
\int_s^t\int_{\mathbb R^d}
\u_r\,dx\,dr
\right)^{\frac1{2q}},
\end{align*}
the diffusion is estimated by
Lemma~\ref{limit:lem:time-integrability}, $p=1+\frac{2}{d(m-1)}$.
All constants are independent of $h,R$ and of the admissible $\psi$ (since its norm is less than $1$).

Let $R\to\infty$. For $0<t-s\le1$, take
$h=|t-s|^{(1-1/p)/2}$ and obtain
\[
d_{BL}(u_t^{\varepsilon,\delta},u_s^{\varepsilon,\delta})
\le C\left(|t-s|^{1/(2q)}
+|t-s|^{(1-1/p)/2}\right).
\]
Since $\frac{1-\frac{1}{p}}{2}=\frac{1}{d(m-1)+2}$, this is the desired estimate.
For larger time differences, enlarge $C$ and use $d_{BL}\le2$.
Theorem~\ref{limit:thm:Ascoli} now gives uniform narrow compactness.

The estimate passes to each $\varepsilon$-limit $u^\delta$.
Applying the same theorem to these curves gives the assertion
for the $\delta$-limit.
This argument for the nonlocal family uses no previously constructed
limit and therefore justifies its earlier use in
Theorem~\ref{limit:thm:epsilon}.
\end{proof}

\begin{theorem}\label{limit:thm:martingale-existence}
Assume $(V),(W),(\mu_0)$ and \eqref{joint}, there
exists a weak solution of \eqref{SDE} whose marginal densities
form a weak solution of \eqref{granular medium equation} with
initial datum $u_0$.
\end{theorem}
\begin{proof}
Let $u$ be the solution obtained in
Theorem~\ref{limit:thm:delta}, with the narrowly continuous
representative identified above.
The limiting estimates imply
\[
\int_0^T\int_{\R^d}
\bigl[u^m+u|\nabla V+\nabla W*u|\bigr]\,dx\,dt<\infty.
\]
Choose a Borel representative of $u$ and set
\[
Lf(t,x)=
-(\nabla V(x)+(\nabla W*u_t)(x))\cdot\nabla f(t,x)
+u_t(x)^{m-1}\Delta f(t,x).
\]
Thus $b=-(\nabla V+\nabla W*u)$ and $a=2u^{m-1}I_d$.
The matrix $a$ is nonnegative, and the displayed integrability
is precisely the condition required by
Theorem~\ref{thm:prelim-superposition}.
Since $u$ solves the corresponding Fokker--Planck equation,
the theorem gives a probability law on $C([0,T];\R^d)$ solving
the martingale problem and having marginals $u_t(x)\,dx$.
The representation of continuous local martingales by stochastic
integrals, on an enlarged probability space if necessary, gives
\[
\begin{aligned}
dY_t&=-(\nabla V(Y_t)+(\nabla W*u_t)(Y_t))\,dt
+\sqrt{2u_t(Y_t)^{m-1}}\,dB_t,\\
\operatorname{Law}(Y_t)(dx)&=u_t(x)\,dx.
\end{aligned}
\]
This proves the assertion.
\end{proof}

% The appendix follows after all main-text sections.

\subsection*{Acknowledgements}
I am grateful to Jan Peszek and Błażej Miasojedow for many helpful suggestions.

\subsection* {Funding}
The author was supported by the National Science Center (Poland), project 2021/43/B/ST1/02851.

\appendix

\section{Extension of the class of test functions}
\label{app:test-functions}

We justify the tests used in Lemma~\ref{limit:lem:entropy-moments}.
The regularity supplied by Corollary~\ref{chaos:cor:nonlocal-existence}
concerns $\sqrt{u^{\varepsilon,\delta}}$ and does not, by itself,
give $u^{\varepsilon,\delta}\in L^2(0,T;H^1(\R^d))$.
We use bounded truncations of the logarithm. A spatial regularization
with the compactly supported mollifier already used in the article
will justify the chain rule; no Gaussian convolution is used.

\begin{lemma}\label{app:lem:entropy-test}
For fixed $\varepsilon,\delta>0$, let $u^{\varepsilon,\delta}$ be
the solution chosen in Corollary~\ref{chaos:cor:nonlocal-existence}.
Then the entropy identity \eqref{final entropy} and the
second-moment identity \eqref{second moment} hold for every
$t\in[0,T]$. In particular, the tests $\log u^{\varepsilon,\delta}+1$
and $|x|^2$ used in the main text are justified.
\end{lemma}

\begin{proof}
    Set
    \begin{align*}
        F&=[\delta+\varphi_\epsilon*(\varphi_\epsilon*\u)^{m-1}]\nabla\u\\&\quad+\u\nabla[\varphi_\epsilon*(\varphi_\epsilon*\u)^{m-1}]+\u(\nabla V+\nabla W*\u).
    \end{align*}

    Now, using Sobolev regularity of solutions obtained in Lemma \ref{limit:lem:fixed-regularity} we can rewrite the weak formulation of \eqref{nonlocal} as
    \begin{align}
        \int_{\R^d}\u_tf_t\,\text{d}x-\int_{\R^d}\u_0f_0\,\text{d}x-\int_0^t\int_{\R^d}\u\partial_rf\,\text{d}x\,\text{d}r=-\int_0^t\int_{\R^d}F\cdot\nabla f\,\text{d}x\,\text{d}r
    \end{align}
    for any $f\in C_c^\infty(\R^d\times[0,T])$.

    Let us fix $\theta,\eta>0$ and denote $\u_\theta:=\u*\varphi_\theta$. Obviously $\partial_t\u_\theta=\text{div}(\varphi_\theta*F)$ in the sense of distributions. This gives the estimate $||\partial_t\u_{\theta}(t)||_{L^1(\R^d)}\leq||\nabla\varphi_\theta||_{L^1(\R^d)}||F||_{L^1(\R^d)}$ which ensures that $\u_\theta\in W^{1,1}([0,T+1];L^1(\R^d))$ since $F$ is in Bochner space $L^1([0,T+1];L^1(\R^d))$ by lemma \ref{limit:lem:fixed-regularity}.

    We choose a solution on $[0,T+1]$ since we need to regularize test function in time and we need solutions on a larger time interval. It exists with the same quantities finite as in the case of interval $[0,T]$. Let $\chi_R$ be a standard family of cut-off functions which are $1$ in ball of radius $R$ and equal to $0$ outside a ball of radius $2R$ with gradient bounded by $\frac{C}{R}$. Now, we introduce a function \[\psi(t,x):=\varphi_\theta*\big[\chi_R\log(\u_\theta(t)+\eta)\big](x)\]
    For fixed $\theta,\eta, R$ this function is smooth in space with compact support. Moreover, it satisfies
    \begin{align}
        \partial_t\psi=\varphi_\theta*\left[\chi_R\frac{\partial_t\u_\theta}{\u_\theta+\eta}\right],
    \end{align}
    so
    \begin{align}\label{regular in time}
||\partial_t\psi(t)||_\infty\leq\frac{||\varphi_\theta||_\infty}{\eta}||\partial_t\u_\theta(t)||_{L^1(\R^d)}
    \end{align}
    which is finite and it ensures that we can mollify in time by some smooth, positive probability kernel $\rho\in C_c^\infty(0,1) $ and obtain a proper family of test functions to the \eqref{nonlocal}:
    \begin{align}
\rho_h(s)=h^{-1}\rho(s/h),
\qquad
\psi_h(t,x)=\int_0^h\rho_h(s)\psi(t+s,x)\,ds.
    \end{align}
    Estimate \eqref{regular in time} gives us the following standard results for mollified sequences:
\begin{align}
&\sup_{r\in[0,T]}
\bigl(\|\psi_h(r)-\psi(r)\|_\infty
+\|\nabla\psi_h(r)-\nabla\psi(r)\|_\infty\bigr)\to0,\\
&\int_0^T\|\partial_r\psi_h-\partial_r\psi\|_\infty\,dr
\to0.
\end{align}
Weak formulation gives us
\begin{align*}
        \int_{\R^d}\u_t\psi_h(t)\,\text{d}x-\int_{\R^d}\u_0\psi_h(0)\,\text{d}x-\int_0^t\int_{\R^d}\u\partial_r\psi_h\,\text{d}x\,\text{d}r=\\-\int_0^t\int_{\R^d}F\cdot\nabla \psi_h\,\text{d}x\,\text{d}r.
        \end{align*}

Since $F$ is in $L^1$ and $\u$ is probability density at every time marginal we can let $h\to0$ by previous estimates. Transferring convolution from $\psi$ to $\u$ we obtain
\begin{align*}
    \int_{\R^d}\u_\theta(t)\chi_R\log(\u_\theta(t)+\eta)\,\text{d}x-\int_{\R^d}\u_\theta(0)\chi_R\log(\u_\theta(0)+\eta)\,\text{d}x\\=\int_0^t\int_{\R^d}\chi_R\frac{\u_\theta}{\u_\theta+\eta}\partial_r\u_\theta\,\text{d}x\,\text{d}r-\int_0^t\int_{\R^d}\chi_R(\varphi_\theta*F)\cdot\frac{\nabla\u_\theta}{\u_\theta+\eta}\,\text{d}x\,\text{d}r\\
    \quad-\int_0^t\int_{\R^d}\log(\u_\theta+\eta)(\varphi_\theta*F)\cdot\nabla\chi_R\,\text{d}x\,\text{d}r
\end{align*}

Let us denote by
\begin{align*}
    \Phi_\eta(r):=(r+\eta)\log(r+\eta)-r-\eta\log\eta.
\end{align*}

Simple calculations allow us to rewrite the above equality us
\begin{align*}
    \int_{\R^d}\chi_R\big[\Phi_\eta(\u_\theta(t))-\Phi_\eta(\u_\theta(0))\big]\,\text{d}x=-\int_0^t\int_{\R^d}\chi_R(\varphi_\theta*F)\cdot\frac{\nabla\u_\theta}{\u_\theta+\eta}\,\text{d}x\,\text{d}r\\
    \quad-\int_0^t\int_{\R^d}\log(\u_\theta+\eta)(\varphi_\theta*F)\cdot\nabla\chi_R\,\text{d}x\,\text{d}r.
\end{align*}
We remove the spatial cut-off. For fixed $\theta,\eta$ logarithm is bounded, hence we can estimate the last term by 
\[
\left|
\int_0^t\int_{\R^d}
\log(\u_\theta+\eta)(\varphi_\theta*F)\cdot\nabla\chi_R
\,\text{d}x\,\text{d}r\right|
\le
\frac{C_{\theta,\eta}}R\|F\|_{L^1}
\to0.
\]
In the remaining terms we converge by Lebesgue Dominated Convergence Theorem. We have
\[
\int_{\R^d}\Phi_\eta(\u_\theta(t))\,\text{d}x
-\int_{\R^d}\Phi_\eta(\u_\theta(0))\,\text{d}x
=
-\int_0^t\int_{\R^d}
(\varphi_\theta*F)\cdot\frac{\nabla \u_\theta}{\u_\theta+\eta}\,\text{d}x\,\text{d}r.
\]
Now we let $\eta\to0$. To do so we use the following estimates:
\[
\begin{aligned}
|(\varphi_\theta*F)(x)|^2
&=
\left|
\int_{\mathbb R^d}
\sqrt{\varphi_\theta(x-y)}\,
\frac{F(y)}{\sqrt{\u(y)}}\,
\sqrt{\varphi_\theta(x-y)\u(y)}\,\text{d}y
\right|^2
\\
&\le
\left(
\int_{\mathbb R^d}
\varphi_\theta(x-y)\frac{|F(y)|^2}{\u(y)}\,dy
\right)
\left(
\int_{\mathbb R^d}
\varphi_\theta(x-y)\u(y)\,\text{d}y
\right)
\\
&=
\left[\varphi_\theta*
\left(\frac{|F|^2}{\u}\right)\right](x)\,\u_\theta(x).
\end{aligned}
\]
Dividing by $\u_\theta$ we get
\begin{align}\label{L2 estimates 1}
    \frac{|(\varphi_\theta*F)(x)|^2}{\u_\theta}\leq\left[\varphi_\theta*\left(\frac{|F|^2}{\u}\right)\right]
\end{align}
and similarly
\begin{align}\label{L2 estimates 2}
\frac{|\nabla\u_\theta|^2}{\u_\theta}\leq\left[\varphi_\theta*\left(\frac{|\nabla\u|^2}{\u}\right)\right].
\end{align}
Both quantities on the right side are integrable in time and space since $2q$-moments of $\u_t$ are bounded for every $t\in[0,T+1]$, fact that $\u$ is probability density for almost every time $t$ and corollary \ref{chaos:cor:nonlocal-existence}. Now, we can estimate using Young's inequality the right-hand side of our equation:
\begin{align*}
(\varphi_\theta*F)\cdot\frac{\nabla \u_\theta}{\u_\theta+\eta}=\frac{\u_\theta}{\u_\theta+\eta}\frac{\varphi_\theta*F}{\sqrt{\u_\theta}}\frac{\nabla\u_\theta}{\sqrt{\u_\theta}}\leq\frac{|(\varphi_\theta*F)(x)|^2}{2\u_\theta}+\frac{|\nabla\u_\theta|^2}{2\u_\theta}.
\end{align*}
This is integrable in time and space by \eqref{L2 estimates 1} and \eqref{L2 estimates 2}. Hence once again by Lebesgue Dominated Convergence Theorem we can pass to the limit on the right-hand side with $\eta\to0$.

On the left-hand side we use
\[
\Phi_\eta(r)\longrightarrow r\log r-r,
\qquad
|\Phi_\eta(r)|\le r|\log r|+Cr
\quad(0<\eta\le1).
\]
This, boundedness of $\u_\theta$ and finite second moment (which ensures that $\u_\theta\log\u_\theta$ is integrable) allow us to converge by Lebesgue Dominated Convergence Theorem. We obtained
\[
H(u_\theta(t))-H(u_\theta(0))
=
-\int_0^t\int
\frac{\varphi_\theta*F}{\sqrt{u_\theta}}
\cdot\frac{\nabla u_\theta}{\sqrt{u_\theta}}
\]
where $H$ is entropy defined in section 1.

The only thing remaining is convergence with $\theta\to0$. We will show it proving strong convergence in $L^2([0,T]\times\R^d)$
\[
\frac{\varphi_\theta*F}{\sqrt{\u_\theta}}\to\frac{F}{\sqrt{\u}},\,\,\,\text{and}\,\,\,\frac{\nabla\u_\theta}{\sqrt{\u_\theta}}\to\frac{\nabla\u}{\sqrt{\u}}.
\]
We will write rigorous proof for the first convergence, the second one is analogous. We know that we have convergence almost everywhere and the estimate \eqref{L2 estimates 1}. Together with  Fatou's lemma, this yields
\begin{align}\label{Fatou}
    4\left|\left|\frac{F}{\sqrt{\u}}\right|\right|_{L^2([0,T]\times\R^d)}^2\leq\liminf_{\theta\to0}\left|\left|\frac{\varphi_\theta*F}{\sqrt{\u_\theta}}+\frac{F}{\sqrt{\u}}\right|\right|^2_{L^2([0,T]\times\R^d)}.
\end{align}
Using Parallelogram law and inequality
\[
\begin{aligned}
\left\|\frac{\varphi_\theta*F}{\sqrt{\u_\theta}}\right\|_{L^2([0,T]\times\R^d)}^2
&\le
\int_0^T\!\int_{\mathbb R^d}
\varphi_\theta*\left(\frac{|F|^2}{\u}\right)\,\text{d}x\,\text{d}s\\
&=
\int_0^T\!\int_{\mathbb R^d}\frac{|F|^2}{\u}\,\text{d}x\,\text{d}s
=
\left\|\frac{F}{\sqrt \u}\right\|_{L^2([0,T]\times\R^d)}^2.
\end{aligned}
\]
give us (where index $2$ denotes $L^2$ norm in $[0,T]\times\R^d$)
\[
\begin{aligned}
0\le{}&
\left\|
\frac{\varphi_\theta*F}{\sqrt{u_\theta}}
-\frac{F}{\sqrt u}
\right\|_2^2\\
={}&
2\left\|\frac{\varphi_\theta*F}{\sqrt{u_\theta}}\right\|_2^2
+2\left\|\frac{F}{\sqrt u}\right\|_2^2
-\left\|
\frac{\varphi_\theta*F}{\sqrt{u_\theta}}
+\frac{F}{\sqrt u}
\right\|_2^2\\
\le{}&
4\left\|\frac{F}{\sqrt u}\right\|_2^2
-\left\|
\frac{\varphi_\theta*F}{\sqrt{u_\theta}}
+\frac{F}{\sqrt u}
\right\|_2^2.
\end{aligned}
\]
Taking $\limsup$ and using \eqref{Fatou} we obtain
\[
\begin{aligned}
0
&\le
\limsup_{\theta\downarrow0}
\left\|
\frac{\varphi_\theta*F}{\sqrt{\u_\theta}}
-\frac{F}{\sqrt \u}
\right\|_2^2\\
&\le
4\left\|\frac{F}{\sqrt \u}\right\|_2^2
-
\liminf_{\theta\downarrow0}
\left\|
\frac{\varphi_\theta*F}{\sqrt{\u_\theta}}
+\frac{F}{\sqrt \u}
\right\|_2^2
\le0.
\end{aligned}
\]
Hence we have desired strong convergence.

We now pass to the limit in the entropy terms. Observe that \[
r\log r+|x|^2r+e^{-1-|x|^2}\geq0.
\]
This with almost everywhere convergence of $\u_\theta\to\u$ allow us to use Fatou lemma
\[
\begin{aligned}
&H(\u)+\int_{\mathbb R^d}|x|^2\u\,\text{d}x
+\int_{\mathbb R^d}e^{-1-|x|^2}\,\text{d}x\\
&\quad\le
\liminf_{\theta\downarrow0}
\left[
H(\u_\theta)+\int_{\mathbb R^d}|x|^2\u_\theta\,\text{d}x
+\int_{\mathbb R^d}e^{-1-|x|^2}\,dx
\right].
\end{aligned}
\]
This with convergence of second moments:
\[
\begin{aligned}
\int_{\mathbb R^d}|x|^2\u_\theta(x)\,\text{d}x
&=\iint_{\mathbb R^d\times\mathbb R^d}
|y+\theta z|^2\u(y)\varphi(z)\,\text{d}y\,\text{d}z\\
&=\int_{\mathbb R^d}|y|^2\u(y)\,\text{d}y
+\theta^2\int_{\mathbb R^d}|z|^2\varphi(z)\,\text{d}z\\
&\longrightarrow\int_{\mathbb R^d}|y|^2\u(y)\,\text{d}y
\end{aligned}
\]
give us inequality
\begin{align}\label{entropy}
    H(\u(t))\leq\liminf_{\theta\to0}H(\u_\theta(t)).
\end{align}
Jensen's inequality gives us 
\begin{align*}
\u_\theta(x)\log\u_\theta(x)\leq\int_{\R^d}\varphi_\theta(y)\u(x-y)\log\u(x-y).
\end{align*}
Integrating with respect to $x$ gives us
\begin{align*}
    H(\u_\theta(t))\leq H(\u(t)).
    \end{align*}
This and \eqref{entropy} give us
\begin{align*}
    H(\u(t))\leq\liminf_{\theta\to0}H(\u_\theta(t))\leq\limsup_{\theta\to0} H(\u_\theta(t))\leq H(\u(t)).
\end{align*}
This yields the desired equality
\begin{align}\label{final entropy}
    H(\u(t))-H(\u(0))
=
-\int_0^t\int_{\mathbb R^d}\frac{F\cdot\nabla \u}{\u}\,\text{d}x\,\text{d}r.
\end{align}

Testing by $|x|^2$ is simpler. We multiply it by $\chi_R$, test by it and use boundedness of all integrals. Lebesgue Dominated Convergence Theorem allows us to obtain
\begin{equation}\label{second moment}
\begin{aligned}
\int |x|^2\u_t&=\int |x|^2\u_0
-2\int_0^t\int x\cdot\nabla V(x)\,\u_r(x)\,dx\,dr\\
&-\int_0^t\iint
(x-y)\cdot\nabla W(x-y)\,\u_r(x)\u_r(y)\,dx\,dy\,dr\\
&+2d\delta t
+2d\int_0^t\int(\varphi_\varepsilon*\u_r)^m\,dx\,dr.
\end{aligned}    
\end{equation}

\end{proof}

\end{document}